# BLOW-UP FOR THE STOCHASTIC NONLINEAR SCHRÖDINGER EQUATION WITH MULTIPLICATIVE NOISE

By Anne de Bouard and Arnaud Debussche

*Université Paris–Sud and ENS de Cachan*

We study the influence of a multiplicative Gaussian noise, white in time and correlated in space, on the blow-up phenomenon in the supercritical nonlinear Schrödinger equation. We prove that any sufficiently regular and localized deterministic initial data gives rise to a solution which blows up in arbitrarily small time with a positive probability.

**1. Introduction.** The understanding of the influence of a noise on the propagation of waves is a very important problem. Although the propagation is often described by deterministic models, in many circumstances randomness should be taken into account. It can change drastically the qualitative behavior and result in new properties.

This is a very vast subject. Propagation can be described by many different models and randomness can take several different forms. In this article, we consider the nonlinear Schrödinger equation as the deterministic model. It describes the propagation of waves in media with both nonlinear and dispersive responses. It is used in many areas of physics, for example, hydrodynamics, plasma physics, nonlinear optics, molecular biology, and so on.

It is well known that this equation has localized solutions called solitary waves or, sometimes, solitons. When the nonlinearity is not too strong—or subcritical, these are particularly robust and propagate without changing form. However, for stronger nonlinearity, these are unstable and the instability results in the collapse—or blow-up—of the localized wave.

Here, we are particularly interested in the influence of a noise acting as a potential on this behavior. Such noise has been considered in [14]; there the

---

Received October 2003; revised March 2004.

*AMS 2000 subject classifications.* Primary 35Q55, 60H15; secondary 76B35, 60H30, 60J60.

*Key words and phrases.* Nonlinear Schrödinger equations, stochastic partial differential equations, white noise, blow-up, variance identity, support theorem.

paths of the noise are smooth functions and the nonlinearity is subcritical. In the case of a white noise, which is considered here, this type of model has been introduced in the context of crystals (see [1, 2] and also [18, 21] for other models). It is expected that such a noise has a strong influence on the solutions which blow-up. It may delay or even prevent the formation of a singularity. In [12], some numerical simulations tend to show that this is the case for a very irregular noise: for a space-time white noise. However, in the supercritical case and for a noise which is correlated in space but nondegenerate, it has been observed that, on the contrary, *any solution* seems to blow-up in a finite time. Recall that in the deterministic case, only a restricted class of solutions blow-up. Our aim is to prove rigorously such a behavior.

It is mathematically very difficult to consider space-time white noises; this is due to the lack of smoothing effect in the Schrödinger equation. Thus, we restrict our attention to the study of correlated noises.

The case of an additive noise has been considered in [8, 9] and it has been proved that *for any initial data*, blow-up occurs in the sense that, for *arbitrary* $t > 0$, the probability that the solution blows up before the time $t$ is strictly positive (see Remark 4.4). Thus, the noise strongly influences this blow-up phenomenon. This result is in perfect agreement with the numerical simulations. The argument is based on three ingredients. First, we generalize the deterministic argument to prove that blow-up occurs for some initial data. This is based on a stochastic version of the variance identity (see [20, 22]). Then, we use the fact that the nonlinear Schrödinger equation is controlable by a forcing term. Thus, any initial data can be transformed into a state which yields a singular solution. Finally, since the noise is nondegenerate and the solution depends continuously on the path of the noise, we can argue that, with positive probability, the noise will be close to the control so that blow-up will happen afterward.

In the multiplicative case, we can again generalize the variance identity and prove that a restricted class of initial data evolves into a singular solution. An initial data with sufficiently negative energy is in this class. This is done in Sections 3 and 4, after we have set some notation and recalled some preliminary results in Section 2. However, it is not known if the nonlinear Schrödinger equation is controlable by a potential and it is well known that the solution does not depend continuously on the noise. Thus, it is not straightforward to generalize the argument used in the additive case. In Section 5 we address this problem and show that, using more subtle arguments, we can prove the same result as in the additive case, when the nonlinearity is supercritical. First, for any initial data, we construct a deterministic potential such that, when it is used as a control in the nonlinear Schrödinger equation, the energy of the solution becomes arbitrarily negative. The problem is then to replace the continuity argument. In fact, we prove that the controlled solution is in the support of the diffusion associated to the stochastic nonlinear equation, which is exactly what we need to prove our result. The techniques used to obtain this last step are inspired from [19] and we refer the reader to this article for a nice and simple presentation of support theorems, as well as for references on this subject. Many ingredients used in [19] are similar to the ones used to prove the convergence of a numerical scheme. Thus, we have been led to use the ideas developed in [11].

**2. Notation.** We consider a complete probability space $(\Omega, \mathcal{F}, \mathbb{P})$, endowed with a filtration $(\mathcal{F}_t)_{t \geq 0}$, and a sequence $(\beta_k)_{k \in \mathbb{N}}$ of independent real-valued Brownian motions on $\mathbb{R}^+$ associated to the filtration $(\mathcal{F}_t)_{t \geq 0}$. We study NLS equations with multiplicative noise. The noise is real valued. To define it rigorously, we use the space $L^2(\mathbb{R}^d; \mathbb{R})$ of real-valued square integrable functions on $\mathbb{R}^d$ and a Hilbertian basis $(e_k)_{k \in \mathbb{N}}$ of this space. Then, given a bounded linear operator $\phi$ on this space, we define the process

$$W(t, x, \omega) = \sum_{k=0}^{\infty} \beta_k(t, \omega) \phi e_k(x), \qquad t \geq 0, \ x \in \mathbb{R}^d, \ \omega \in \Omega.$$

It is a Wiener process on $L^2(\mathbb{R}^d; \mathbb{R})$, with covariance operator $\phi \phi^*$. In all that follows, $\phi \phi^*$ will be assumed to be a finite trace operator in $L^2(\mathbb{R}^d; \mathbb{R})$; this easily implies, for a fixed $t$, the convergence of the series above in $L^2(\Omega; L^2(\mathbb{R}^d; \mathbb{R}))$, and almost surely in $L^2(\mathbb{R}^d; \mathbb{R})$. The time derivative of $W$ models a noise which is delta correlated in time. If $\phi$ is defined through a kernel $\mathcal{K}$, which means that for any square integrable function $u$,

$$\phi u(x) = \int_{\mathbb{R}^d} \mathcal{K}(x, y) u(y) \, dy,$$

then the correlation function of the noise is formally given by

$$\mathbb{E} \left( \frac{\partial W}{\partial t}(t, x) \frac{\partial W}{\partial t}(s, y) \right) = c(x, y) \delta_{t-s},$$

with

$$c(x, y) = \int_{\mathbb{R}^d} \mathcal{K}(x, z) \mathcal{K}(y, z) \, dz.$$

Let $\sigma > 0$, we consider the following stochastic NLS equation, introduced in [1] and [2],

$$(2.1) \qquad idu - (\Delta u + |u|^{2\sigma} u) \, dt = u \circ dW,$$

where $\circ$ stands for a Stratonovich product in the right-hand side of (2.1). We will use the equivalent Itô form of this equation. Define, for $x \in \mathbb{R}^d$, the function

$$F_\phi(x) = \sum_{k=0}^{\infty} (\phi e_k(x))^2,$$

which does not depend on the basis $(e_k)_{k \in \mathbb{N}}$; this equivalent Itô equation may be written as

$$(2.2) \qquad i\,du - (\Delta u + |u|^{2\sigma} u)\,dt = u\,dW - \frac{i}{2} u F_\phi\,dt$$

(see [7] for details).

A local existence result for (2.2) has been proved in [10]. Before recalling it, we need to set the functional framework.

For $p \geq 1$, $L^p(\mathbb{R}^d)$ is the classical Lebesgue space (of complex valued functions), and the inner product in the real Hilbert space $L^2(\mathbb{R}^d)$ is denoted by $(\cdot, \cdot)$, that is,

$$(u, v) = \operatorname{Re} \int_{\mathbb{R}^d} u(x) \bar{v}(x)\,dx$$

for $u, v \in L^2(\mathbb{R}^d)$. The norm in $L^p(\mathbb{R}^d)$ is denoted by $|\cdot|_{L^p}$ and, for $p = 2$, by $|\cdot|$.

We define the usual space $H^s(\mathbb{R}^d)$ of tempered distributions $v \in \mathcal{S}'(\mathbb{R}^d)$, whose Fourier transform $\hat{v}$ satisfies $(1 + |\xi|^2)^{s/2} \hat{v} \in L^2(\mathbb{R}^d)$; $H^s(\mathbb{R}^d, \mathbb{R})$ is the subspace of $H^s(\mathbb{R}^d)$ consisting of real-valued functions. The norm in $H^s(\mathbb{R}^d)$ is denoted by $|\cdot|_{H^s}$. For $p \in \mathbb{N}$, $W^{1,p}(\mathbb{R}^d)$ is the space of functions in $L^p(\mathbb{R}^d)$ whose first order derivatives are in $L^p(\mathbb{R}^d)$ and its norm is denoted by $|\cdot|_{W^{1,p}}$.

If $I$ is an interval of $\mathbb{R}$, $E$ is a Banach space, and $1 \leq r \leq +\infty$, then $L^r(I; E)$ is the space of strongly Lebesgue measurable functions $v$ from $I$ into $E$ such that the function $t \mapsto |v(t)|_E$ is in $L^r(I)$. We define similarly the spaces $L^r(\Omega; E)$.

Given two separable Hilbert spaces $H$ and $\widetilde{H}$, we denote by $\mathcal{L}_2(H, \widetilde{H})$ the space of Hilbert–Schmidt operators $\Phi$ from $H$ into $\widetilde{H}$, endowed with the norm

$$|\Phi|^2_{\mathcal{L}_2(H, \widetilde{H})} = \operatorname{tr} \Phi^* \Phi = \sum_{k \in \mathbb{N}} |\Phi e_k|^2_{\widetilde{H}},$$

where $(e_k)_{k \in \mathbb{N}}$ is any orthonormal basis of $H$. When $H = L^2(\mathbb{R}^d; \mathbb{R})$ and $\widetilde{H} = H^s(\mathbb{R}^d; \mathbb{R})$, then $\mathcal{L}_2(H, \widetilde{H})$ is simply denoted by $\mathcal{L}_2^{0,s}$. Given a Banach space $B$, we will also consider bounded linear operators from $L^2(\mathbb{R}^d)$ into $B$, and in order to replace the notion of Hilbert–Schmidt operators, we use in this case the notion of $\gamma$-radonifying operators (see, e.g., [3, 4]). We denote by $R(L^2, B)$ the space of $\gamma$-radonifying operators from $L^2(\mathbb{R}^d; \mathbb{R})$ into $B$. We recall that if $B$ is a Hilbert space, then $R(L^2, B)$ coincides with $\mathcal{L}_2(L^2(\mathbb{R}^d; \mathbb{R}), B)$.

It is classical that if $\psi(\cdot)$ is a predictable random process defined on $[0, T]$, with values in the space of continuous operators from $L^2(\mathbb{R}^d, \mathbb{R})$ into an

Hilbert space $K$, and is such that

$$\int_0^T |\psi(s)\phi|^2_{\mathcal{L}_2(L^2(\mathbb{R}^d,\mathbb{R}),K)}\,ds < \infty \qquad \text{a.s.,}$$

then the stochastic integral

$$\int_0^T \psi(s)\,dW(s) = \sum_{\ell \in \mathbb{N}} \int_0^T \psi(s)\phi e_\ell\,d\beta_\ell(s)$$

is a well-defined $K$-valued random variable.

We now recall the existence result of [10]. For technical reasons, we restrict our attention to the case $d \leq 3$. However, if the assumptions are suitably modified, most of the results can be extended to larger spatial dimensions.

We recall that a pair $(r,p)$ of positive numbers is called an admissible pair if $r \geq 2$ and $\frac{2}{r} = d(\frac{1}{2} - \frac{1}{p})$.

**Theorem 2.1.** *Assume that $0 < \sigma$ if $d = 1$ or $2$ and $0 < \sigma < 2$ if $d = 3$. Let $\phi \in \mathcal{L}_2^{0,1} \cap R(L^2, W^{1,\alpha}(\mathbb{R}^d))$ with $\alpha > 2d$. Then there exists an admissible pair $(r,p)$ such that for any $\mathcal{F}_0$-measurable $u_0$ with values in $H^1(\mathbb{R}^d)$, there exist a stopping time $\tau^*(u_0,\omega)$ and a unique solution to (2.2) starting from $u_0$ which is almost surely in $C([0,\tau]; H^1(\mathbb{R}^d)) \cap L^r(0,\tau; W^{1,p}(\mathbb{R}^d))$ for any $\tau < \tau^*(u_0)$. Moreover, we have almost surely,*

$$\tau^*(u_0,\omega) = +\infty \quad or \quad \limsup_{t \nearrow \tau^*(u_0,\omega)} |u(t)|_{H^1(\mathbb{R}^d)} = +\infty.$$

In this article we want to study blow-up phenomena. As in the deterministic case, we need to work in the context of $H^1$ solutions. This is why we use the result of [10]. In [7], $L^2$ solutions are considered, requiring weaker assumptions on $\phi$; however, there, a stronger restriction is imposed on $\sigma$ and it is not expected that blow-up occurs in that case. Indeed, it is shown in [7] that the solutions are global in time. In fact, we will require that $\sigma$ does not fulfill the assumption of [7]. We will work in the so-called critical and supercritical cases:

$$\sigma \geq \frac{2}{d}.$$

The critical case corresponds to $\sigma = \frac{2}{d}$ and most of our results will not be true in that case. Note that if $\sigma$ satisfies the assumptions of Theorem 2.1, by the Sobolev embedding theorem, we know that $H^1(\mathbb{R}^d) \subset L^{2\sigma+2}(\mathbb{R}^d)$ with a continuous embedding.

In fact, we need stronger assumptions on the initial data. As in most of the deterministic blow-up results, we require some spatial localization. For $\eta \geq 0$, we introduce the spaces

$$\Sigma^\eta = \{v \in H^\eta(\mathbb{R}^d) : |x|^\eta v \in L^2(\mathbb{R}^d)\}$$

endowed with the norm $|\cdot|_{\Sigma^\eta}$:

$$|v|^2_{\Sigma^\eta} = |v|^2_{H^\eta} + ||x|^\eta v|^2.$$

When $\eta = 1$, we set $\Sigma^1 = \Sigma$.

In all of the article, $c$ or $c(\cdot, \ldots, \cdot)$ is a constant which may change from one line to another and depends only on its arguments.

**3. The stochastic variance identity.** In this section we derive an identity on the evolution of the variance, or virial, of the solutions of the nonlinear Schrödinger equation (2.2). This is a generalization of the well-known corresponding formula in the deterministic case.

The variance (this quantity should not be confused with the probabilistic variance) is defined by

$$V(v) = \int_{\mathbb{R}^d} |x|^2 |v(x)|^2 \, dx, \qquad v \in \Sigma.$$

Its evolution is described in terms of the mass,

$$M(v) = \int_{\mathbb{R}^d} |v(x)|^2 \, dx = |v|^2, \qquad v \in L^2(\mathbb{R}^d),$$

the energy,

$$H(v) = \frac{1}{2} \int_{\mathbb{R}^d} |\nabla v(x)|^2 \, dx - \frac{1}{2\sigma + 2} \int_{\mathbb{R}^d} |v(x)|^{2\sigma+2} \, dx$$

$$= \frac{1}{2} |\nabla v|^2 - \frac{1}{2\sigma + 2} |v|^{2\sigma+2}_{L^{2\sigma+2}}, \qquad v \in H^1(\mathbb{R}^d),$$

and the momentum,

$$G(v) = \mathrm{Im} \int_{\mathbb{R}^d} v(x) \, x \cdot \nabla \bar{v}(x) \, dx, \qquad v \in \Sigma.$$

The first result describes the evolution of $M$ and $H$ and is proved in [10]. We introduce the function $f^1_\phi$ depending on the covariance operator $\phi$ and such that, for any orthonormal basis $(e_\ell)_{\ell \in \mathbb{N}}$ of $L^2(\mathbb{R}^d; \mathbb{R})$,

$$f^1_\phi(x) = \sum_{\ell \in \mathbb{N}} |\nabla \phi e_\ell(x)|^2, \qquad x \in \mathbb{R}^d.$$

Note that this function does not depend on the chosen basis $(e_\ell)_{\ell \in \mathbb{N}}$. Indeed, since $\phi$ is assumed to be Hilbert–Schmidt from $L^2(\mathbb{R}^d)$ into $H^1(\mathbb{R}^d)$, it is associated to a kernel $\mathcal{K} \in L^2(\mathbb{R}^d \times \mathbb{R}^d)$ such that $\nabla_x \mathcal{K} \in L^2(\mathbb{R}^d \times \mathbb{R}^d)$; it can be seen that $f^1_\phi(x) = |\nabla_x \mathcal{K}(x, \cdot)|^2_{L^2(\mathbb{R}^d)}$.

PROPOSITION 3.1. *Let $u_0$, $\sigma$ and $\phi$ be as in Theorem* 2.1. *For any stopping time $\tau$ such that $\tau < \tau^*(u_0)$ a.s., we have*

$$(3.1) \qquad M(u(\tau)) = M(u_0) \qquad a.s.$$

*and*

$$
\begin{aligned}
(3.2) \quad H(u(\tau)) = H(u_0) - \mathrm{Im} \sum_{\ell \in \mathbb{N}} \int_0^\tau \int_{\mathbb{R}^d} \bar{u}(s,x) \nabla u(s,x) \cdot \nabla(\phi e_\ell)(x) \, dx \, d\beta_\ell(s) \\
+ \tfrac{1}{2} \int_0^\tau \int_{\mathbb{R}^d} |u(s,x)|^2 f_\phi^1(x) \, dx \, ds \qquad a.s.,
\end{aligned}
$$

*where $u$ is the solution of* (2.2) *given by Theorem* 2.1 *with $u(0) = u_0$.*

The different terms in (3.2) make sense. To see this we first observe that, since $\phi \in \mathcal{L}_2^{0,1}$, the function $f_\phi^1$ is an $L^1$ function and since $\phi \in R(L^2, W^{1,\alpha}(\mathbb{R}^d))$, we also have $f_\phi^1 \in L^{\alpha/2}(\mathbb{R}^d)$ (see [3, 4, 7]). Thus, thanks to Hölder's inequality, $f_\phi^1 \in L^q(\mathbb{R}^d)$ for any $q \in [1, \alpha/2]$. We then introduce the stopping time:

$$\tau_k = \inf\{s \in [0,t], \ |u(s)|_{H^1} \geq k\}.$$

For any $k \in \mathbb{N}$, we have

$$
\mathbb{E}\left( \left| \sum_{\ell \in \mathbb{N}} \int_0^{\tau \wedge \tau_k} \int_{\mathbb{R}^d} \bar{u}(s,x) \nabla u(s,x) \cdot \nabla(\phi e_\ell)(x) \, dx \, d\beta_\ell(s) \right|^2 \right)
$$

$$
= \mathbb{E}\left( \sum_{\ell \in \mathbb{N}} \int_0^{\tau \wedge \tau_k} \left| \int_{\mathbb{R}^d} \bar{u}(s,x) \nabla u(s,x) \cdot \nabla(\phi e_\ell)(x) \, dx \right|^2 ds \right)
$$

$$
\leq \mathbb{E}\left( \int_0^{\tau \wedge \tau_k} \left( \int_{\mathbb{R}^d} |\bar{u}(s,x)| |\nabla u(s,x)| (f_\phi^1(x))^{1/2} \, dx \right)^2 ds \right),
$$

thanks to Minkowski's inequality in the last step. Using Sobolev's embedding, it is now easy to see that this is a finite quantity. Since, for almost every $\omega \in \Omega$, we have $\tau \leq \tau_k$ for some $k$, it follows that the stochastic term in (3.2) is a well-defined stochastic integral.

We now investigate the evolution of $V$ and $G$.

PROPOSITION 3.2. *Let the assumptions of Proposition* 3.1 *hold and assume, moreover, that $u_0 \in \Sigma$ a.s., then for any stopping time $\tau$ such that $\tau < \tau^*(u_0)$ a.s., we have*

$$
\begin{aligned}
(3.3) \quad G(u(\tau)) = G(u_0) + 4 \int_0^\tau H(u(s)) \, ds + \frac{2 - \sigma d}{\sigma + 1} \int_0^\tau |u(s)|_{L^{2\sigma+2}}^{2\sigma+2} \, ds \\
+ \sum_{k \in \mathbb{N}} \int_0^\tau \int_{\mathbb{R}^d} |u(s,x)|^2 x \cdot \nabla(\phi e_k)(x) \, dx \, d\beta_k(s)
\end{aligned}
$$

*and*

$$(3.4) \qquad V(u(\tau)) = V(u_0) + 4 \int_0^\tau G(u(s)) \, ds.$$

REMARK 3.3. The proof of this result is postponed to Section 6.1. Note that, formally, it is the result of the same computation as in the deterministic case and that it is easier to perform the computation with the Stratonovich equation since in that case the standard calculus can be used. However, in the proof we use the Itô form which is better suited for rigorous justifications.

A similar result was obtained in the case of an additive noise in [8]. The formula was more complicated there. In the case of a multiplicative noise, more terms cancel in the computation. This strongly uses the fact that the noise is real valued and reflects the fact that a Stratonovich multiplicative noise corresponds to a noisy potential and is a more realistic representation of physical phenomena in the present context.

REMARK 3.4. We will actually prove that, for each $t \geq 0$ and $k \in \mathbb{N}$, there exists a constant $c(V(u_0), t, k)$ such that $V(u(t)) \leq c(V(u_0), t, k)$ provided $t < \tau^*(u_0)$ and $|u(s)|_{H^1} \leq k$ for $s \in [0, t]$ [see (6.2)].

Using this and similar arguments as above, it can be seen that the stochastic term in (3.3) is well defined.

COROLLARY 3.5. *Under the same assumptions as in Proposition* 3.2, *we have*

$$
\begin{aligned}
V(u(\tau)) = \; & V(u_0) + 4G(u_0)\tau + 8H(u_0)\tau^2 \\
& + 4 \frac{2 - \sigma d}{\sigma + 1} \int_0^\tau \int_0^s |u(s_1)|^{2\sigma+2}_{L^{2\sigma+2}} \, ds_1 \, ds \\
& + 8 \int_0^\tau \int_0^s \int_0^{s_1} \int_{\mathbb{R}^d} |u(s_2, x)|^2 f_\phi^1(x) \, dx \, ds_2 \, ds_1 \, ds \\
& + 4 \sum_{k \in \mathbb{N}} \int_0^\tau \int_0^s \int_{\mathbb{R}^d} |u(s_1, x)|^2 x \cdot \nabla(\phi e_k)(x) \, dx \, d\beta_k(s_1) \, ds \\
& - 16 \operatorname{Im} \sum_{k \in \mathbb{N}} \int_0^\tau \int_0^s \int_0^{s_1} \int_{\mathbb{R}^d} \bar{u}(s_2, x) \nabla u(s_2, x) \\
& \qquad\qquad\qquad\qquad \times \nabla(\phi e_k)(x) \, dx \, d\beta_k(s_2) \, ds_1 \, ds.
\end{aligned}
$$

(3.5)

**4. Blow-up for initial data with negative energy.** The aim of this section is to generalize the well-known deterministic result stating that, for the critical or supercritical NLS equation, an initial data with negative energy yields a solution which forms a singularity in finite time. (See [20], Chapter 5, where more sophisticated results can also be found.)

THEOREM 4.1. *Let $u_0$, $\sigma$ and $\phi$ satisfy the assumptions of Theorem 2.1. Assume also that $\sigma \geq \frac{2}{d}$, $u_0 \in L^2(\Omega; \Sigma) \cap L^{2\sigma+2}(\Omega; L^{2\sigma+2}(\mathbb{R}^d))$, $f_\phi^1$ is a bounded function and, for some $\bar{t} > 0$,*

$$(4.1) \quad \mathbb{E}(V(u_0)) + 4\mathbb{E}(G(u_0))\bar{t} + 8\mathbb{E}(H(u_0))\bar{t}^2 + \frac{4}{3}\bar{t}^3 m_\phi \mathbb{E}(M(u_0)) < 0,$$

*where $m_\phi = |f_\phi^1|_{L^\infty}$; then*

$$\mathbb{P}(\tau^*(u_0) \leq \bar{t}) > 0.$$

REMARK 4.2. Clearly, if the energy of $u_0$ is a.s. negative, then the conclusion of Theorem 4.1 holds for some $\bar{t} > 0$ provided the noise is not too strong, that is, provided $m_\phi$ is small enough.

REMARK 4.3. By Sobolev's embedding theorem, the condition that $f_\phi^1$ is bounded is fulfilled if $\phi$ is a $\gamma$-radonifying operator from $L^2(\mathbb{R}^d; \mathbb{R})$ into $W^{s,p}(\mathbb{R}^d; \mathbb{R})$ with $(s-1)p > d$.

REMARK 4.4. In [8], in the case of an additive noise, this result was proved under a severe assumption on $\sigma$ and the proof was rather technical. Under the assumptions made here, only a weaker result was proved. The argument used below can, in fact, be easily adapted to the additive case and the severe assumption on $\sigma$ can be removed.

PROOF OF THEOREM 4.1. Assume that the conclusion of Theorem 4.1 does not hold; then $\bar{t} < \tau^*(u_0)$ a.s. and we can take $\tau = \bar{t}$ as a stopping time in Corollary 3.5. In that case (3.5) can be simplified into

$$
\begin{aligned}
(4.2) \quad V(u(\bar{t})) &= V(u_0) + 4G(u_0)\bar{t} + 8H(u_0)\bar{t}^2 \\
&\quad + 4\frac{2-\sigma d}{\sigma+1} \int_0^{\bar{t}} (\bar{t}-s)|u(s)|_{L^{2\sigma+2}}^{2\sigma+2} \, ds \\
&\quad + 4\int_0^{\bar{t}} (\bar{t}-s)^2 \int_{\mathbb{R}^d} |u(s,x)|^2 f_\phi^1 \, dx \, ds \\
&\quad + 4\sum_{\ell \in \mathbb{N}} \int_0^{\bar{t}} (\bar{t}-s) \int_{\mathbb{R}^d} |u(s,x)|^2 x \cdot \nabla(\phi e_\ell)(x) \, dx \, d\beta_\ell(s) \\
&\quad - 8\,\mathrm{Im} \sum_{\ell \in \mathbb{N}} \int_0^{\bar{t}} (\bar{t}-s)^2 \int_{\mathbb{R}^d} \bar{u}(s,x) \nabla u(s,x) \cdot \nabla(\phi e_\ell)(x) \, dx \, d\beta_\ell(s).
\end{aligned}
$$

Let us set for $t, r \geq 0$,

$$V(t,r) = V(u_0) + 4G(u_0)t + 8H(u_0)t^2$$

$$+ 4\frac{2 - \sigma d}{\sigma + 1} \int_0^t (t - s)|u(s)|^{2\sigma + 2}_{L^{2\sigma + 2}} \, ds$$

$$+ 4\int_0^t (t - s)^2 \int_{\mathbb{R}^d} |u(s, x)|^2 f_\phi^1 \, dx \, ds$$

$$+ 4\sum_{\ell \in \mathbb{N}} \int_0^r (t - s) \int_{\mathbb{R}^d} |u(s, x)|^2 x \cdot \nabla(\phi e_\ell)(x) \, dx \, d\beta_\ell(s)$$

$$- 8\, \mathrm{Im} \sum_{\ell \in \mathbb{N}} \int_0^r (t - s)^2 \int_{\mathbb{R}^d} \bar{u}(s, x) \nabla u(s, x) \cdot \nabla(\phi e_\ell)(x) \, dx \, d\beta_\ell(s),$$

clearly, $V(u(\bar{t})) = V(\bar{t}, \bar{t})$. Define the stopping time

$$\tau_k = \inf\{s \in [0, \bar{t}], \ |u(s)|_{H^1} \geq k\},$$

for $k \in \mathbb{N}$. Since the mass of $u$ is constant, we have, assuming $t \leq \bar{t}$,

$$\mathbb{E}\left(\sum_{\ell \in \mathbb{N}} \int_0^{\tau_k} \left|(t - s)\int_{\mathbb{R}^d} |u(s, x)|^2 x \cdot \nabla(\phi e_\ell)(x) \, dx\right|^2 ds\right)$$

$$\leq m_\phi \mathbb{E}\left(\int_0^{\tau_k} (t - s)^2 M(u_0) V(u(s)) \, ds\right)$$

$$\leq \tfrac{1}{3} m_\phi \bar{t}^3 \sup_{s \in [0, \bar{t}]} \mathbb{E}(V(u(s \wedge \tau_k)) M(u_0))$$

since $m_\phi = |f_\phi^1|_{L^\infty}$. By Remark 3.4, we deduce that this quantity is finite. It follows that, when $r = \tau_k$, the first stochastic integral above is square integrable and its expectation vanishes. Similarly,

$$\mathbb{E}\left(\sum_{\ell \in \mathbb{N}} \int_0^{\tau_k} \left|(t - s)^2 \int_{\mathbb{R}^d} \bar{u}(s, x) \nabla u(s, x) \cdot \nabla(\phi e_\ell)(x) \, dx\right|^2 ds\right)$$

$$\leq \tfrac{1}{5} m_\phi \bar{t}^5 k^2 \mathbb{E}(M(u_0))$$

and, when $r = \tau_k$, the expectation of the last term also vanishes. Moreover, using again the conservation of the mass,

$$\int_0^t (t - s)^2 \int_{\mathbb{R}^d} |u(s, x)|^2 f_\phi^1(x) \, dx \, ds \leq \tfrac{1}{3} t^3 m_\phi M(u_0)$$

and, recalling that $\sigma \geq \frac{2}{d}$, we get that, for any $k \in \mathbb{N}$,

$$\mathbb{E}(V(t, \tau_k)) \leq \mathbb{E}(V(u_0)) + 4\mathbb{E}(G(u_0))t + 8\mathbb{E}(H(u_0))t^2 + \tfrac{4}{3}t^3 m_\phi \mathbb{E}(M(u_0)).$$

We now choose $t = \bar{t}$. By assumption, $\tau_k \to \bar{t}$ a.s. as $k$ goes to infinity. Thus, thanks to Fatou's lemma, we have

$$\mathbb{E}(V(u(\bar{t}))) \leq \mathbb{E}(V(u_0)) + 4\mathbb{E}(G(u_0))\bar{t} + 8\mathbb{E}(H(u_0))\bar{t}^2 + \tfrac{4}{3}\bar{t}^3 m_\phi \mathbb{E}(M(u_0)).$$

This last inequality contradicts (4.1), since $V(u(\bar{t}))$ is nonnegative. Thus, the result is proved.  □

The result stated in Theorem 4.1 is still valid if the expectation over $\Omega$ in (4.1) is replaced by the expectation over an $\mathcal{F}_0$-measurable subset of $\Omega$. More precisely, defining, for $\bar{M}$ and $\bar{H}$ positive constants,

(4.3)  $\mathcal{V}_{\bar{M},\bar{H}} = \{v \in \Sigma, \ V(v) < \bar{M}, \ G(v) < \bar{M}, \ |v|_{L^2}^2 < \bar{M}, \ H(v) < -\bar{H}\},$

we can prove the following.

COROLLARY 4.5.   *Let $u_0$, $\sigma$ and $\phi$ satisfy the assumptions of Theorem 2.1. Assume also that $\sigma \geq \frac{2}{d}$, $u_0 \in \Sigma$ a.s. and $f_\phi^1$ is a bounded function. Then for any $\bar{M} > 0$ and $\bar{t} > 0$, there is a constant $\bar{H}(\bar{t}, \bar{M}) > 0$ such that*

$$\mathbb{P}(u_0 \in \mathcal{V}_{\bar{M},\bar{H}}) > 0 \quad \implies \quad \mathbb{P}(\tau^*(u_0) \leq t) > 0.$$

PROOF.   This is simply proved by taking $\bar{H}$ large enough so that

$$\bar{M} + 4\bar{t}\bar{M} - 8\bar{t}^2\bar{H} + \tfrac{4}{3}\bar{t}^3 m_\phi \bar{M} < 0$$

and applying Theorem 4.1 with $u_0$ replaced by $u_0 \mathbb{1}_{\Omega_0}$, where

$$\Omega_0 = \{\omega \in \Omega, \ u_0(\cdot, \omega) \in \mathcal{V}_{\bar{M},\bar{H}}\}.$$

Note, indeed, that the solution of (2.2) satisfies $u(t, \omega) = 0$ for any $t \geq 0$ if $u(0, \omega) = 0$.  □

**5. Blow-up for any initial data in the supercritical case.**   In this section we assume that the nonlinearity is supercritical, $\sigma > \frac{2}{d}$. In that case we are able to strengthen considerably the result of Theorem 4.1.

5.1. *Main result.*   The main result of this article is the following.

THEOREM 5.1.   *Assume that $\frac{2}{d} < \sigma$ if $d = 1, 2$ and $\frac{2}{3} < \sigma < 2$ if $d = 3$. Let $\phi \in \mathcal{L}_2^{0,2}$ be such that $\ker \phi^* = \{0\}$. Then for any $u_0 \in \Sigma^2$ with $u_0 \neq 0$ and $t > 0$, we have*

$$\mathbb{P}(\tau^*(u_0) < t) > 0,$$

*where $\tau^*(u_0)$ is the existence time of the solution of (2.2) with initial data $u_0$ provided by Theorem 2.1.*

The proof of this result is given in Section 5.4. The idea is the same as in [8] where the case of an additive noise is treated. We choose $t_1$, $t_2 > 0$ and consider a control problem where the noise is replaced by a control and show that, for any initial data, there exists a control such that the solution $U$ of the controlled nonlinear Schrödinger equation at time $t_1$ satisfies

$$(5.1) \quad V(U(t_1)) + 4G(U(t_1))t_2 + 8H(U(t_1))t_2^2 + \frac{4}{3}t_2^3 m_\phi M(U(t_1)) < 0.$$

Then, using the nondegeneracy of the noise, that is, $\ker \phi^* = \{0\}$, we know that the noise will be close to the control on the time interval $[0, t_1]$ with positive probability. Intuitively, this means that at time $t_1$ the solution of (2.2) will be close to $U(t_1)$ on a set of positive probability and that we will have $\mathbb{P}(u(t_1) \in \mathcal{V}_{\bar{M}, \bar{H}}) > 0$ for some $\bar{M}$ and $\bar{H}$ as in Corollary 4.5—the initial time being $t_1$ instead of 0 and $\bar{t} = t_2$. Thus, for the solution starting at 0, blow-up will occur with positive probability before the time $t_1 + t_2$. This is exactly the result since $t_1$, $t_2$ are arbitrary.

In the case of a multiplicative noise considered here, this is much more difficult to justify rigorously.

First, the control is a potential and it is not known whether the nonlinear Schrödinger equation is controlable by a potential. We prove in Section 5.2 that, for any initial data, it is possible to construct a potential such that the solution at arbitrary time $t_1$ has a very negative energy so that (5.1) holds. This construction only works in the supercritical case. It is not difficult to see that in the critical case, $\sigma = \frac{2}{d}$, such a result cannot be true.

The second difficulty is that the Itô or Stratonovich products are not continuous with respect to the paths of the noise. Then, to replace the continuity argument used in [8], we use more sophisticated probabilistic tools inspired by the characterization of the support of a diffusion (see [19]). We prove that, if the noise is nondegerate, the solution of the controlled problem is in the support of the diffusion associated with (2.2), which implies that with positive probability the solution will be close to the solution of the controlled problem. This is the aim of Section 5.3. Finally, we end the proof in Section 5.4. Technical results are proved in Section 6.

The assumption on the smoothness of $u_0$ is quite strong here since in the deterministic literature it is common to consider initial data in $\Sigma$. Our result can probably be extended to such data, however, this would considerably complicate the proof in Section 6.2 and, for clarity, we have decided to restrict our attention to this smaller class of initial data.

5.2. *Construction of a potential leading to blow-up.* We first show that for a supercritical nonlinearity it is possible to build a control which acts as a potential in the nonlinear Schrödinger equation and has the property that, in a finite time, the solution has a very negative energy. Therefore, if the control is then switched off, and replaced by the noise, it follows from Corollary 4.5 that the solution will rapidly blow-up.

PROPOSITION 5.2. *Let $u_0 \in \Sigma$, $u_0 \neq 0$, $T_1 > 0$. Let $\bar{M} = \max(M(u_0) + 1; V(u_0) + 4|G(u_0)| + 1)$. Then, for any $\bar{H} \geq 0$ there exists $T_2 \leq T_1$ and a potential $f \in L^s(0, T_2; W^{1,p}(\mathbb{R}^d))$, for some $s > 1$ and some $p$ with $1 \leq p < 1 + \frac{1}{\sigma}$, such that the solution of*

$$(5.2) \qquad i\frac{dU}{dt} - (\Delta U + |U|^{2\sigma}U + fU) = 0,$$
$$U(0) = u_0$$

*exists on $[0, T_2]$ and satisfies $U(T_2) \in \mathcal{V}_{\bar{M}, \bar{H}}$, where $\mathcal{V}_{\bar{M}, \bar{H}}$ is defined by (4.3).*

*Moreover, if $u_0 \in H^2(\mathbb{R}^d)$, then $f \in C([0, T_2]; H^2(\mathbb{R}^d; \mathbb{R}))$.*

PROOF. We can assume that $T_1 \leq 1$. Let $\frac{2}{d} < \tilde{\sigma} < \sigma$. We then consider the following supercritical nonlinear Schrödinger equation:

$$(5.3) \qquad i\frac{dU}{dt} - (\Delta U + \lambda|U|^{2\tilde{\sigma}}U) = 0,$$
$$U(0) = u_0.$$

The constant $\lambda \geq 0$ is taken such that the energy corresponding to equation (5.3) satisfies

$$H_{\tilde{\sigma}, \lambda}(u_0) = \int_{\mathbb{R}^d} \frac{1}{2}|\nabla u_0(x)|^2 - \frac{\lambda}{2\tilde{\sigma}+2}|u_0(x)|^{2\tilde{\sigma}+2}\, dx$$
$$< -\frac{1}{8T_1^2}(V(u_0) + 4T_1|G(u_0)|).$$

Since $\tilde{\sigma} < \sigma$, and $\sigma < 2$ if $d = 3$, it follows from the deterministic theory (see [5]) that there exists a solution $U$ in $C([0, t^*]; \Sigma) \cap L^r(0, t^*; W^{1,p}(\mathbb{R}^d))$ for some $t^* > 0$ and $(r, p)$ admissible with $p = 2\tilde{\sigma} + 2$. Moreover, by the assumption above,

$$V(u_0) + 4T_1 G(u_0) + 8T_1^2 H_{\tilde{\sigma}, \lambda}(u_0) < 0$$

and, since $\frac{2}{d} < \tilde{\sigma}$, this implies that $U$ blows up at some time $T^* < T_1$. This is classical ([15, 20, 22]) and follows from the deterministic variance identity which states that if $U$ exists on the time interval $[0, T_1]$, then

$$(5.4) \qquad V(U(t)) \leq V(u_0) + 4tG(u_0) + 8t^2 H_{\tilde{\sigma}, \lambda}(u_0); \qquad t \in [0, T_1].$$

But since $V(U(t))$ is nonnegative, this is clearly impossible. In fact, we have

$$|U(t)|_{H^1} \to +\infty \quad \text{and} \quad |U(t)|_{L^{2\tilde{\sigma}+2}} \to +\infty \qquad \text{when } t \to T^*.$$

By the conservation of the $L^2$ norm and Hölder's inequality, we have

$$|U(t)|_{L^{2\tilde{\sigma}+2}}^{2\tilde{\sigma}+2} \leq \frac{2\tilde{\sigma}+2}{2\lambda(2\sigma+2)}|U(t)|_{L^{2\sigma+2}}^{2\sigma+2} + c(\lambda, \sigma, \tilde{\sigma})|u_0|_{L^2}^2;$$

it follows

$$H(U(t)) = \frac{\lambda}{2\tilde{\sigma} + 2}|U(t)|_{L^{2\tilde{\sigma}+2}}^{2\tilde{\sigma}+2} - \frac{1}{2\sigma + 2}|U(t)|_{L^{2\sigma+2}}^{2\sigma+2} + H_{\tilde{\sigma},\lambda}(U(t))$$

$$\leq -\frac{\lambda}{2\tilde{\sigma} + 2}|U(t)|_{L^{2\tilde{\sigma}+2}}^{2\tilde{\sigma}+2} + \frac{2\lambda}{2\tilde{\sigma} + 2}c(\lambda,\sigma,\tilde{\sigma})|u_0|_{L^2}^2 + H_{\tilde{\sigma},\lambda}(U(t)).$$

Recalling that $H_{\tilde{\sigma},\lambda}(U(t))$ does not depend on $t$, we deduce that

$$H(U(t)) \to -\infty \qquad \text{when } t \to T^*.$$

Thus, there exists $T_2 \leq T_1$ such that

$$H(U(T_2)) < -\bar{H}.$$

Moreover, by (5.4),

$$V(U(T_2)) \leq V(u_0) + 4T_2 G(u_0) + 8T_2^2 H_{\tilde{\sigma},\lambda}(u_0)$$

$$\leq V(u_0) + 4|G(u_0)| < \bar{M}$$

since $T_2 \leq 1$ and $H_{\tilde{\sigma},\lambda}(u_0) \leq 0$. It is also classical that

$$G(U(T_2)) \leq G(u_0) + 4T_2 H_{\tilde{\sigma},\lambda}(u_0),$$

so that $G(U(T_2)) \leq G(u_0) < \bar{M}$. In addition, $M(u(T_2)) = M(u_0) < \bar{M}$, hence, $U(T_2) \in \mathcal{V}_{\bar{M},\bar{H}}$. It suffices to take $f(t) = \lambda|U(t)|^{2\tilde{\sigma}} - |U(t)|^{2\sigma}$ and $U$ is then the solution of (5.2).

Finally, we observe that if $u_0 \in H^2(\mathbb{R}^d)$, it follows from [17] that $U \in C([0, T_2]; H^2(\mathbb{R}^d))$. Now, if $d = 3$, $f$ is not necessarily in $\in C([0, T_2]; H^2(\mathbb{R}^d))$, but it is easily seen that $f \in C([0, T_2]; W^{1,p}(\mathbb{R}^d))$ for any $p$ with $2 \leq p < 2d/(d-2)$, and we may argue by density as follows: take $f_n \in C([0, T_2]; H^2(\mathbb{R}^d))$ with $f_n \to f$ in $C([0, T_2]; W^{1,p}(\mathbb{R}^d))$, where $p < 2d/(d-2)$ is such that for some $q < 2d/(d-2)$, we have $\frac{1}{q'} = 1 - \frac{1}{q} = \frac{1}{2} + \frac{1}{p}$, that is, $\frac{1}{p} + \frac{1}{q} = \frac{1}{2}$. This is possible thanks to the fact that $d < 4$. Using then the method in [17], it is not difficult to see that the solution $U_n$ of (5.2), with $f$ replaced by $f_n$, exists on $[0, T_2]$ and $U_n \to U$ in $C([0, T_2]; H^1(\mathbb{R}^d))$. In addition, $xU_n \to xU$ in $C([0, T_2]; L^2(\mathbb{R}^d))$, as follows from the same arguments as in [5], Corollary 6.4.4. Hence, for $n$ sufficiently large, $f_n$ fits the conclusion of Proposition 5.2. $\square$

5.3. *The support result.* We now show that, if the noise is nondegenerate, the solution $U$ given by Proposition 5.2 is in the support of the diffusion associated to (2.2). More precisely, under some additional assumptions on the initial data, we prove that the support of the law of $u(T_2)$ contains a neighborhood in $\Sigma$ of $U(T_2)$. In fact, following [19] and slightly generalizing our arguments, we could characterize exactly the support of the diffusion $u$. However, this is not the aim of this paper and, for clarity, we simply state and prove the result we need.

PROPOSITION 5.3. *Suppose that the assumptions of Theorem 2.1 hold,* $u_0 \in \Sigma^2$, $u_0 \neq 0$, $\phi \in \mathcal{L}_2^{0,2}$ *and the noise is nondegenerate:* $\ker \phi^* = \{0\}$. *Let* $T_1 > 0$, $\bar{H} > 0$ *arbitrary, and* $T_2 \leq T_1$, $U$ *be given by Proposition 5.2; let* $u$ *be the solution of (2.2) given by Theorem 2.1; then for arbitrary neighborhood* $\mathcal{V}$ *of* $U(T_2)$ *in* $\Sigma$, *we have*

$$\mathbb{P}(\tau^*(u_0) > T_2 \text{ and } u(T_2) \in \mathcal{V}) > 0.$$

PROOF. We adapt the method used to characterize the support of a diffusion and follow the same line as in [19].

We choose $T > T_2$. Let $(e_k)_{k \in \mathbb{N}}$ be a complete orthonormal system in $L^2(\mathbb{R}^d; \mathbb{R})$ and $(\beta_k)_{k \in \mathbb{N}}$ be a sequence of independent Brownian motions such that $W = \sum_{k \in \mathbb{N}} \beta_k \phi e_k$. We define $W_c = \sum_{k \in \mathbb{N}} \beta_k e_k$ so that $W_c$ is a cylindrical Wiener process on $L^2(\mathbb{R}^d; \mathbb{R})$.

Since $\ker \phi^* = \{0\}$, we know that the range of $\phi$ is dense in $L^2(\mathbb{R}^d, \mathbb{R})$ and thus in $H^2(\mathbb{R}^d, \mathbb{R})$. It follows that, for any $n \in \mathbb{N}$, there exists $g_n \in C([0, T_2]; L^2(\mathbb{R}^d, \mathbb{R}))$ such that

$$(5.5) \qquad |f - \phi g_n|_{C([0,T_2];H^2(\mathbb{R}^d))} \leq \frac{1}{n},$$

where $f$ is the potential given by Proposition 5.2. We set

$$f_n(t) = \phi g_n(t), \qquad t \leq T_2, \qquad f_n(t) = 0, \qquad T_2 < t \leq T.$$

We also assume that we have chosen $(e_k)_{k \in \mathbb{N}}$ in such a way that $e_k \in H^2(\mathbb{R}^d)$, $k \in \mathbb{N}$. We denote by $P_n$ the orthogonal projector onto $Sp(e_0, \ldots, e_n)$ and set $\Delta t = \frac{1}{n}$. We construct a piecewise constant approximation of the noise in the following way:

$$\dot{W}_{c,n}(t) = \frac{P_n W_c(k\Delta t) - P_n W_c((k-1)\Delta t)}{\Delta t}, \qquad t \in [k\Delta t, (k+1)\Delta t[,$$

and

$$\dot{W}_n = \phi \dot{W}_{c,n}$$

[we set $W_c(-\Delta t) = 0$]. Then, we consider the following equation, for $n \in \mathbb{N}$:

$$(5.6) \quad i \, du^n - (\Delta u^n + |u^n|^{2\sigma} u^n + f_n u^n) \, dt = u^n \, dW - \frac{i}{2} u^n F_\phi \, dt - u^n \dot{W}_n \, dt.$$

By Girsanov's theorem (see [6], Theorem 10.14),

$$W_{c,n}(t) = W_c(t) - \int_0^t (\dot{W}_{c,n}(s) - g_n(s)) \, ds$$

is a cylindrical Wiener process for the probability measure

$$d\mathbb{P}_n = D_n \, d\mathbb{P},$$

where

$$D_n = \exp\left[\int_0^T (\dot{W}_{c,n}(s) - g_n(s), dW_c(s)) - \frac{1}{2}\int_0^T |\dot{W}_{c,n}(s) - g_n(s)|^2\, ds\right].$$

Clearly, $u^n$ is the solution of

$$(5.7) \qquad i\,du^n - (\Delta u^n + |u^n|^{2\sigma}u^n)\,dt = u^n\,dW_n - \frac{i}{2}u^n F_\phi\,dt,$$

with $W_n = \phi W_{c,n}$. Since the law of the solution of the stochastic nonlinear Schrödinger equation does not depend on the probability space nor on the Wiener process, we deduce that the law of $u^n$ is the same as the law of $u$, the solution of (2.2). Therefore, for any Borelian set $\mathcal{V}$,

$$\mathbb{P}(\tau^*(u_0) > T_2 \text{ and } u(T_2) \in \mathcal{V}) = \mathbb{P}_n(\tau_n^*(u_0) > T_2 \text{ and } u^n(T_2) \in \mathcal{V})$$

$$= \int_\Omega \mathbb{1}_{\{\tau_n^*(u_0) > T_2 \text{ and } u^n(T_2) \in \mathcal{V}\}} D_n\, d\mathbb{P},$$

where we have denoted by $\tau_n^*(u_0)$ the existence time for the solution of (5.7) with initial data $u_0$. It follows

$$\mathbb{P}(\tau^*(u_0) > T_2 \text{ and } u(T_2) \in \mathcal{V}) > 0$$

$$\text{if and only if } \mathbb{P}(\tau_n^*(u_0) > T_2 \text{ and } u^n(T_2) \in \mathcal{V}) > 0.$$

We prove in Section 6.2 that

$$(5.8) \qquad\qquad \lim_{n\to\infty} \mathbb{P}(\tau_n^*(u_0) \le T_2) = 0$$

and

$$(5.9) \quad \forall \epsilon > 0 \qquad \lim_{n\to\infty} \mathbb{P}\left(\tau_n^*(u_0) > T_2 \text{ and } \sup_{t\in[0,T_2]} |u^n(t) - U(t)|_\Sigma > \epsilon\right) = 0.$$

Therefore, for any $\mathcal{V}$ neighborhood in $\Sigma$ of $U(T_2)$,

$$\mathbb{P}(\tau_n^*(u_0) > T_2 \text{ and } u^n(T_2) \in \mathcal{V}) \to 1 \qquad \text{when } n \to \infty,$$

which clearly implies the result. $\quad\square$

5.4. *Proof of Theorem* 5.1. It is now easy to complete the argument sketched in Section 5.1. Indeed, for any $t > 0$ and $u_0 \in \Sigma^2$, we choose $t_1$, $t_2 > 0$ such that $t_1 + t_2 \le t$. We take $\bar{M}$, as in Proposition 5.2, and choose $\bar{H}$ such that Corollary 4.5 is satisfied with $\bar{M}$, and $\bar{t} = t_2$.

Then we take $T_1 = t_1$ in Proposition 5.2 and obtain the control $f$ and the solution $U$ on the interval $[0, T_2]$ with $T_2 \le t_1$. It follows from the statement of Proposition 5.2 that $U(T_2) \in \mathcal{V}_{\bar{M}, \bar{H}}$.

We then use Proposition 5.3, with $\mathcal{V} = \mathcal{V}_{\bar{M}, \bar{H}}$; setting

$$\Omega_{T_2} = \{\omega \in \Omega \,|\, \tau^*(u_0) > T_2 \text{ and } u(T_2) \in \mathcal{V}\},$$

which is clearly an $\mathcal{F}_{T_2}$-measurable set, and which has positive probability according to Proposition 5.3, we deduce by Corollary 4.5, that the solution starting at time $0$ with initial data $u(T_2)$ will blow-up before $t_2$ with a positive probability. But, if we denote this latter solution by $\tilde{u}(t)$, then it is clear, by translating the Brownian motions, that the solution $u$ of (2.2) satisifies for $t \geq T_2 : u(T_2 + \cdot)$ and $\tilde{u}(\cdot)$ have the same law. Hence, $u$ blows up with positive probability before $T_2 + t_2$, and the result follows.

## 6. Proofs of technical results.

6.1. *Proof of Proposition* 3.2. The proof is based on a smoothing procedure, as in [10], and on the use of a truncated form of the variance:

$$V_\varepsilon(v) = \int_{\mathbb{R}^d} e^{-\varepsilon|x|^2}|x|^2|v(x)|^2\,dx, \qquad v \in L^2(\mathbb{R}^d).$$

The smoothing step is similar to the one in the proof of Proposition 3.1, which can be found in [10]. Thus, we omit it, keeping in mind that the computation below can be rigorously justified in this way.

We apply Itô's formula to $V_\varepsilon(u(t))$ (see [6], Theorem 4.17 for the form of Itô formula used here):

$$\begin{aligned}
dV_\varepsilon(u(t)) = {}&(V_\varepsilon'(u(t)), -i\Delta u(t))\,dt + (V_\varepsilon'(u(t)), -i|u(t)|^{2\sigma}u(t))\,dt \\
&- \tfrac{1}{2}(V_\varepsilon'(u(t)), u(t)F_\phi)\,dt - (V_\varepsilon'(u(t)), iu(t)\,dW(t)) \\
&+ \tfrac{1}{2}\mathrm{Tr}(V_\varepsilon''(u(t))(iu(t)\phi)(iu(t)\phi)^*)\,dt.
\end{aligned}$$

Note that, for $v \in L^2(\mathbb{R}^d)$,

$$(V_\varepsilon'(u(t)), v) = 2\mathrm{Re}\left(\int_{\mathbb{R}^d} e^{-\varepsilon|x|^2}|x|^2 u(t,x)\bar{v}(x)\,dx\right).$$

Hence, this term cancels if $v$ is a real-valued function multiplied by $iu(t)$. This is the case of the second and fourth terms above.

We have

$$(V_\varepsilon'(u(t)), -i\Delta u(t)) = 4\mathrm{Im}\int_{\mathbb{R}^d} e^{-\varepsilon|x|^2}(1-\varepsilon|x|^2)(x \cdot \nabla u(t,x))\bar{u}(t,x)\,dx,$$

as follows after integrating by parts. Moreover,

$$\begin{aligned}
\tfrac{1}{2}\mathrm{Tr}&(V_\varepsilon''(u(t))(iu(t)\phi)(iu(t)\phi)^*) \\
&= \sum_{k\in\mathbb{N}} \int_{\mathbb{R}^d} e^{-\varepsilon|x|^2}|x|^2(iu(t)\phi e_k)(x)(-i\bar{u}(t)\phi e_k)(x)\,dx \\
&= \int_{\mathbb{R}^d} e^{-\varepsilon|x|^2}|x|^2|u(t,x)|^2 F_\phi(x)\,dx
\end{aligned}$$

and

$$\tfrac{1}{2}(V'_\varepsilon(u(t)), u(t)F_\phi)\, dt = \int_{\mathbb{R}^d} e^{-\varepsilon|x|^2}|x|^2|u(t,x)|^2 F_\phi(x)\, dx.$$

It follows that the contributions of these two terms cancel. (This cancellation is natural. It reflects the fact that the noise is real and acts as a potential.) We obtain

$$dV_\varepsilon(u(t)) = 4\mathrm{Im}\int_{\mathbb{R}^d} e^{-\varepsilon|x|^2}(1-\varepsilon|x|^2)(x\cdot\nabla u(t,x))\bar{u}(t,x)\, dx\, dt$$

and for $t < \tau^*(u_0)$,

(6.1)
$$\begin{aligned}
V_\varepsilon(u(t)) &= V_\varepsilon(u_0) \\
&\quad + 4\,\mathrm{Im}\int_0^t\int_{\mathbb{R}^d} e^{-\varepsilon|x|^2}(1-\varepsilon|x|^2)(x\cdot\nabla u(s,x))\bar{u}(s,x)\, dx\, ds.
\end{aligned}$$

In particular, if

$$\tau_k = \inf\{t\in[0,T],\ |u(t)|_{H^1}\geq k\},$$

we have by the Cauchy–Schwarz inequality

$$V_\varepsilon(u(t\wedge\tau_k)) \leq V_\varepsilon(u_0) + 4k\int_0^{t\wedge\tau_k} V_\varepsilon^{1/2}(u(s))\, ds$$

and by Gronwall's lemma

$$V_\varepsilon(u(t\wedge\tau_k)) \leq (4k^2 T + V_\varepsilon(u_0))e^T.$$

Letting $\varepsilon\to 0$, this proves that, for $t<\tau_k$, $u(t)\in\Sigma$ and

(6.2)
$$V(u(t\wedge\tau_k)) \leq (4k^2 T + V(u_0))e^T.$$

Thus, $u$ has trajectories in $L^\infty(0,\tau_k;\Sigma)$ for any $k\in\mathbb{N}$. We now let $k\to+\infty$ and deduce that $u$ has trajectories in $L^\infty(0,\tau;\Sigma)$ for any stopping time $\tau < \tau^*(u_0)$. It is now easy to let $\varepsilon\to 0$ in (6.1) and prove (3.4).

Similarly, we apply Itô's formula to $G(u(t))$ and use the identity

$$\begin{aligned}
(G'(u(t)), &-i(\Delta u(t) + |u(t)|^{2\sigma}u(t))) \\
&= (2ix\cdot\nabla u(t) + diu(t), -i(\Delta u(t) + |u(t)|^{2\sigma}u(t))) \\
&= 2\int_{\mathbb{R}^d}|\nabla u(t,x)|^2\, dx - \frac{\sigma d}{\sigma+1}\int_{\mathbb{R}^d}|u(t,x)|^{2\sigma+2}\, dx \\
&= 4H(u(t)) + \frac{2-\sigma d}{\sigma+1}|u(t)|_{L^{2\sigma+2}}^{2\sigma+2}.
\end{aligned}$$

We also have, after integration by parts,

$$
(G'(u(t)), -iu\, dW(t))
$$
$$
= \operatorname{Re} \sum_{\ell \in \mathbb{N}} \int_{\mathbb{R}^d} (2ix \cdot \nabla u(t,x) + diu(t,x)) i\bar{u}(t,x)(\phi e_\ell)(x)\, dx\, d\beta_\ell(t)
$$
$$
= \sum_{\ell \in \mathbb{N}} \int_{\mathbb{R}^d} |u(t,x)|^2 x \cdot \nabla(\phi e_\ell)(x)\, dx\, d\beta_\ell(t).
$$

Finally, a similar computation shows that

$$
(G'(u(t)), -\tfrac{1}{2} u(t) F_\phi) + \tfrac{1}{2} \operatorname{Tr}(G''(u(t))(iu(t)\phi)(iu(t)\phi)^*) = 0.
$$

Gathering these results and integrating in time yields (3.3).

6.2. *Proof of* (5.8) *and* (5.9). We first introduce a truncated form of (5.6). We choose $3/2 < s_0 < 2$ and a cut-off function $\vartheta \in C_0^\infty(\mathbb{R})$ such that

$$
\begin{aligned}
\vartheta(x) &= 1, & x &\in [0,1], \\
\vartheta(x) &\in [0,1], & x &\geq 0, \\
\vartheta(x) &= 0, & x &\geq 2.
\end{aligned}
$$

Then, for any $R \geq 1$, we set

$$
\theta_R(u) = \vartheta\left(\frac{|u|_{H^{s_0}}^2}{R^2}\right), \qquad u \in H^{s_0}(\mathbb{R}^d).
$$

We then consider the following equation:

$$
\begin{aligned}
(6.3) \quad & idu_R^n - \Delta u_R^n\, dt - \theta_R(u_R^n)(|u_R^n|^{2\sigma} u_R^n + f_n u_R^n)\, dt \\
& = \theta_R(u_R^n) u_R^n\, dW + \Lambda_{R,\phi}(u_R^n)\, dt - \theta_R(u_R^n) u_R^n \dot{W}_n\, dt,
\end{aligned}
$$

with initial data

$$
u_R^n(0) = u_0.
$$

We have denoted by

$$
\Lambda_{R,\phi}(u_R^n) = -\frac{i}{2}\theta_R^2(u_R^n) F_\phi u_R^n - \frac{1}{R^2} \sum_{\ell \in \mathbb{N}} \theta_R(u_R^n) \vartheta'\left(\frac{|u|_{H^{s_0}}^2}{R^2}\right)(u, iu\phi e_\ell)_{H^{s_0}} u\phi e_\ell.
$$

Note that, formally,

$$
\theta_R(u_R^n) u_R^n\, dW + \Lambda_{R,\phi}(u_R^n)\, dt = \theta_R(u_R^n) u_R^n \circ dW.
$$

It is not difficult to see that, for any $R \geq 1$, (6.3) has a unique solution with paths in $C([0,T]; H^2(\mathbb{R}^d))$.

STEP 1.  *For any $R \geq 1$, there exists a constant $C(R, T, \phi)$ such that, for any $n \in \mathbb{N}$,*

$$\mathbb{E}\left(\sup_{t \in [0,T]} |u_R^n(t)|_{H^2}^2\right) \leq C(R, T, \phi).$$

PROOF.   In order to lighten the notation in this step, we omit to explicit the dependence on $R$. All the constants appearing below are allowed to depend on $R$, $T$ or $\phi$, but not on $n$. We set

$$u^n = u_R^n, \qquad \theta = \theta_R, \qquad \Lambda = \Lambda_{R,\phi}.$$

We use the Itô formula with $F(u^n) = |u^n|_{H^2}^2 = |(I - \Delta)u^n|^2$ and we have

$$d|u^n|_{H^2}^2 + 2\mathrm{Re}\left[i \int_{\mathbb{R}^d} \theta(u^n)((I - \Delta)^2 \bar{u}^n)(|u^n|^{2\sigma}u^n + f_n u^n)\, dx\right] dt$$

$$= 2\mathrm{Re}\left(-i \int_{\mathbb{R}^d} \theta(u^n)((I - \Delta)^2 \bar{u}^n) u^n\, dW\, dx\right)$$

$$+ 2\mathrm{Re}\left(-i \int_{\mathbb{R}^d} ((I - \Delta)^2 \bar{u}^n)\Lambda(u^n)\, dx\right) dt$$

$$+ 2\mathrm{Re}\left(i \int_{\mathbb{R}^d} \theta(u^n)((I - \Delta)^2 \bar{u}^n) u^n \dot{W}_n\, dx\right) dt + \theta^2(u^n) \sum_{\ell \in \mathbb{N}} |u^n \phi e_\ell|_{H^2}^2\, dt.$$

We integrate between $0$ and $t$, and estimate each term of the resulting identity. For the second term of the left-hand side, we write

$$\mathbb{E}\left(\sup_{t \in [0,T]} \int_0^t \left| \int_{\mathbb{R}^d} \theta(u^n)((I - \Delta)^2 \bar{u}^n)(|u^n|^{2\sigma}u^n + f_n u^n)\, dx \right| ds\right)$$

$$\leq \mathbb{E}\left(\int_0^T \left| \int_{\mathbb{R}^d} \theta(u^n)(I - \Delta)\bar{u}^n (I - \Delta)(|u^n|^{2\sigma}u^n + f_n u^n)\, dx \right| ds\right)$$

$$\leq c\mathbb{E}\left(\int_0^T \theta(u^n)|u^n|_{H^2}(||u^n|^{2\sigma}u^n|_{H^2} + |f_n u^n|_{H^2})\, ds\right).$$

We have used an integration by parts and Cauchy–Schwarz's inequality. Now, we use the inequality

$$|vw|_{H^2} \leq |v|_{L^\infty}|w|_{H^2} + |w|_{L^\infty}|v|_{H^2}, \qquad v, w \in H^2(\mathbb{R}^d),$$

and the Sobolev embedding $H^{s_0}(\mathbb{R}^d) \subset L^\infty(\mathbb{R}^d)$ to get that

$$(6.4) \qquad |vw|_{H^2} \leq |v|_{H^{s_0}}|w|_{H^2} + |w|_{H^{s_0}}|v|_{H^2}, \qquad v, w \in H^2(\mathbb{R}^d).$$

This gives

$$\theta(u^n)||u^n|^{2\sigma}u^n|_{H^2} \leq cR^{2\sigma}|u^n|_{H^2} \leq c|u^n|_{H^2}$$

for another constant $c$ which depends on $R$. By (5.5), $f_n$ is in $C([0,T]; H^2(\mathbb{R}^d))$ and its norm in this space is bounded independently of $n$. Therefore, by a similar argument as above,

$$|f_n u^n|_{H^2} \le c|u^n|_{H^2}.$$

We deduce

$$\mathbb{E}\bigg( \sup_{t\in[0,T]} \int_0^t \bigg| \int_{\mathbb{R}^d} \theta(u^n)((I-\Delta)^2 \bar{u}^n)(|u^n|^{2\sigma}u^n + f_n u^n)\,dx \bigg|\,ds \bigg)$$

$$\le c\mathbb{E}\bigg( \int_0^T |u^n|_{H^2}^2\,ds \bigg).$$

In the same way, using $\phi \in \mathcal{L}_2^{0,2}$, we have

$$\mathbb{E}\bigg( \sup_{t\in[0,T]} \int_0^t \bigg| \int_{\mathbb{R}^d} ((I-\Delta)^2 \bar{u}^n)\Lambda(u^n)\,dx \bigg|\,ds \bigg) \le c\mathbb{E}\bigg( \int_0^T |u^n|_{H^2}^2\,ds \bigg)$$

and

$$\mathbb{E}\bigg( \sup_{t\in[0,T]} \int_0^t \bigg| \sum_{\ell\in\mathbb{N}} \theta^2(u^n)|u^n\phi e_\ell|_{H^2}^2\,ds \bigg| \bigg) \le c\mathbb{E}\bigg( \int_0^T |u^n|_{H^2}^2\,ds \bigg).$$

Moreover, recalling that the noise is real valued and using a martingale inequality,

$$\mathbb{E}\bigg( \sup_{t\in[0,T]} \bigg| \mathrm{Re} \int_0^t \int_{\mathbb{R}^d} -i\theta(u^n)((I-\Delta)^2 \bar{u}^n)u^n\,dx\,dW \bigg| \bigg)$$

$$= \mathbb{E}\bigg( \sup_{t\in[0,T]} \bigg| \int_0^t \mathrm{Im} \sum_{\ell\in\mathbb{N}} \int_{\mathbb{R}^d} \theta(u^n)(I-\Delta)\bar{u}^n(I-\Delta)(u^n\phi e_\ell)\,dx\,d\beta_\ell \bigg| \bigg)$$

$$= \mathbb{E}\bigg( \sup_{t\in[0,T]} \bigg| \int_0^t \mathrm{Im} \sum_{\ell\in\mathbb{N}} \int_{\mathbb{R}^d} \theta(u^n)(I-\Delta)\bar{u}^n$$

$$\times [u^n\Delta(\phi e_\ell) + 2\nabla u^n \cdot \nabla(\phi e_\ell)]\,dx\,d\beta_\ell \bigg| \bigg)$$

$$\le c\bigg( \mathbb{E}\int_0^T \sum_{\ell\in\mathbb{N}} \bigg( \int_{\mathbb{R}^d} \theta(u^n)(I-\Delta)\bar{u}^n$$

$$\times [u^n\Delta(\phi e_\ell) + 2\nabla u^n \cdot \nabla(\phi e_\ell)]\,dx \bigg)^2\,dt \bigg)^{1/2}$$

$$\le c\bigg( \mathbb{E}\int_0^T \theta(u^n) \sum_{\ell\in\mathbb{N}} (|u^n|_{L^\infty}|u^n|_{H^2}|\phi e_\ell|_{H^2}$$

$$+ |u^n|_{H^2}|\nabla u^n|_{L^3}|\nabla\phi e_\ell|_{L^6})^2\,ds \bigg)^{1/2},$$

by Hölder's inequality. We now use the embeddings $H^{s_0}(\mathbb{R}^d) \subset L^\infty(\mathbb{R}^d)$, $H^{s_0-1}(\mathbb{R}^d) \subset H^{1/2}(\mathbb{R}^d) \subset L^3(\mathbb{R}^d)$ and $H^1(\mathbb{R}^d) \subset L^6(\mathbb{R}^d)$ to obtain

$$\mathbb{E}\Bigg( \sup_{t\in[0,T]} \Big| \mathrm{Re} \int_0^t (-i) \int_{\mathbb{R}^d} \theta(u^n)((I-\Delta)^2 \bar{u}^n) u^n \, dx \, dW \Big| \Bigg)$$

$$\leq c \Bigg( \mathbb{E} \int_0^T \sum_{\ell \in \mathbb{N}} \theta(u^n) |u^n|^2_{H^{s_0}} |u^n|^2_{H^2} |\phi e_\ell|^2_{H^2} \, dt \Bigg)^{1/2}$$

$$\leq c + \tfrac{1}{2}\mathbb{E}\Big( \sup_{t\in[0,T]} |u^n|^2_{H^2} \Big).$$

It remains to estimate the most difficult term which we denote by $A^n(t)$. By similar manipulations as above, we have

$$A^n(t) = \int_0^t \mathrm{Re}\Big( i \int_{\mathbb{R}^d} \theta(u^n)((I-\Delta)^2 \bar{u}^n) u^n \dot{W}_n \, dx \Big) \, ds$$

$$= -\int_0^t \mathrm{Im}\Big( \int_{\mathbb{R}^d} \theta(u^n)(I-\Delta)\bar{u}^n(u^n \Delta \dot{W}_n + 2\nabla u^n \cdot \nabla \dot{W}_n) \, dx \Big) \, ds.$$

We need to introduce a further notation. For $s \geq 0$ and $k \in \mathbb{N}$, we set $s_k = k\Delta t$ and $u^n_S(s, s_{k-1}) = S(s - s_{k-1}) u^n(s_{k-1})$. Here $(S(t))_{t\in\mathbb{R}}$ is the group of operators associated to the linear Schrödinger equation, $S(t) = e^{-it\Delta}$, $t \in \mathbb{R}$. It is well known that, for any $t$, $S(t)$ is an isometry on each Sobolev space $H^s(\mathbb{R}^d)$, $s \in \mathbb{R}$. Now we split $A^n(t)$ as follows:

$$-A^n(t) = \int_0^t \mathrm{Im}\Big( \int_{\mathbb{R}^d} \theta(u^n(s))(I-\Delta)(\bar{u}^n(s) - \bar{u}^n_S(s, s_{k-1}))$$

$$\times (u^n \Delta \dot{W}_n + 2\nabla u^n \cdot \nabla \dot{W}_n) \, dx \Big) \, ds$$

$$+ \int_0^t \mathrm{Im}\Big( \int_{\mathbb{R}^d} \theta(u^n(s))(I-\Delta)\bar{u}^n_S(s, s_{k-1})$$

$$\times (2\nabla(u^n(s) - u^n_S(s, s_{k-1})) \cdot \nabla \dot{W}_n) \, dx \Big) \, ds$$

$$+ \int_0^t \mathrm{Im}\Big( \int_{\mathbb{R}^d} \theta(u^n(s))(I-\Delta)\bar{u}^n_S(s, s_{k-1})$$

$$\times ((u^n(s) - u^n_S(s, s_{k-1})) \Delta \dot{W}_n) \, dx \Big) \, ds$$

$$+ \int_0^t \mathrm{Im}\Big( \int_{\mathbb{R}^d} (\theta(u^n(s)) - \theta(u^n_S(s, s_{k-1})))(I-\Delta)\bar{u}^n_S(s, s_{k-1})$$

$$\times (u^n_S(s, s_{k-1}) \Delta \dot{W}_n + 2\nabla u^n_S(s, s_{k-1}) \cdot \nabla \dot{W}_n) \, dx \Big) \, ds$$

$$+ \int_0^t \text{Im}\left(\int_{\mathbb{R}^d} \theta(u_S^n(s, s_{k-1}))(I - \Delta)\bar{u}_S^n(s, s_{k-1})\right.$$

$$\left. \times (u_S^n(s, s_{k-1})\Delta\dot{W}_n + 2\nabla u_S^n(s, s_{k-1}) \cdot \nabla\dot{W}_n) \, dx\right) ds,$$

and in this expression, $s_{k-1}$ depends on $s$ and is chosen to satisfy $s \in [s_{k-1}, s_k[$. We use similar arguments as above to deduce

$$|A^n(t)| \leq c \int_0^t (|u^n(s) - u_S^n(s, s_{k-1})|_{H^2}|\dot{W}_n|_{H^2}$$

$$+ |u_S^n(s, s_{k-1})|_{H^2}|u^n(s) - u_S^n(s, s_{k-1})|_{H^{s_0}}|\dot{W}_n|_{H^2}$$

$$+ |u^n(s) - u_S^n(s, s_{k-1})|_{H^{s_0}}|u_S^n(s, s_{k-1})|_{H^{s_0}}$$

$$\times |u_S^n(s, s_{k-1})|_{H^2}|\dot{W}_n|_{H^2}) \, ds$$

$$+ \left|\int_0^t \text{Im}\left(\int_{\mathbb{R}^d} \theta(u_S^n(s, s_{k-1}))(I - \Delta)\bar{u}_S^n(s, s_{k-1})\right.\right.$$

$$\left.\left. \times (u_S^n(s, s_{k-1})\Delta\dot{W}_n + 2\nabla u_S^n(s, s_{k-1}) \cdot \nabla\dot{W}_n) \, dx\right) ds\right|$$

$$= |A_1^n(t)| + |A_2^n(t)|.$$

The second term $A_2^n$ is a martingale and can be treated as the stochastic integral above. We obtain

$$\mathbb{E}\left(\sup_{t \in [0,T]} A_2^n(t)\right)$$

$$\leq c\left(\mathbb{E}\int_0^T \sum_{\ell \in \mathbb{N}} \theta^2(u_S^n(s, s_{k-1}))|u_S^n(s, s_{k-1})|_{H^{s_0}}^2|u_S^n(s, s_{k-1})|_{H^2}^2|\phi e_\ell|_{H^2}^2 \, dt\right)^{1/2}$$

$$\leq c + \tfrac{1}{4}\mathbb{E}\left(\sup_{t \in [0,T]} |u^n|_{H^2}^2\right).$$

The term $A_1^n$ is more difficult to treat since $\mathbb{E}(|\dot{W}_n|_{H^2}^m)$ is of order $(\Delta t)^{-m/2}$ for any $m \in \mathbb{N}$. This factor is compensated by $u^n(s) - u_S^n(s, s_{k-1})$, which we now estimate. We write, using the integral form of (6.3) for $u^n$, starting at $t = s_{k-1}$,

$$u^n(s) - u_S^n(s, s_{k-1})$$

$$= i\int_{s_{k-1}}^s S(s - \tau)[\theta(u^n(\tau))|u^n(\tau)|^{2\sigma}u^n(\tau) + f_n(\tau)u^n(\tau) + \Lambda(u^n(\tau))] \, d\tau$$

$$+ i\int_{s_{k-1}}^s S(s - \tau)\theta(u^n(\tau))u^n(\tau) \, dW(\tau)$$

$$- i \int_{s_{k-1}}^{s} S(s-\tau)\theta(u^n(\tau))u^n(\tau)\dot{W}_n \, d\tau.$$

Since $H^{s_0}(\mathbb{R}^d)$ is an algebra, we obtain

$$|u^n(s) - u_S^n(s, s_{k-1})|_{H^{s_0}}$$
$$\leq c(s-s_{k-1}) + \left| \int_{s_{k-1}}^{s} S(s-\tau)\theta(u^n(\tau))u^n(\tau) \, dW(\tau) \right|_{H^{s_0}}$$
$$+ c(s-s_{k-1})|\dot{W}_n((k-1)\Delta t)|_{H^{s_0}}$$

and, for any $m \in \mathbb{N}$,

$$(6.5) \qquad \mathbb{E}(|u^n(s) - u_S^n(s, s_{k-1})|_{H^{s_0}}^m) \leq c(m)(s-s_{k-1})^{m/2}.$$

Moreover, by (6.4),

$$|u^n(s) - u_S^n(s, s_{k-1})|_{H^2}$$
$$\leq c(s-s_{k-1}) \sup_{[0,s]} |u^n|_{H^2} + \left| \int_{s_{k-1}}^{s} S(s-\tau)\theta(u^n(\tau))u^n(\tau) \, dW(\tau) \right|_{H^2}$$
$$+ c(s-s_{k-1}) \sup_{[0,s]} |u^n|_{H^2} |\dot{W}_n((k-1)\Delta t)|_{H^2}$$

and, by Cauchy–Schwarz's inequality,

$$\mathbb{E}(|u^n(s) - u_S^n(s, s_{k-1})|_{H^2} |\dot{W}_n(s)|_{H^2}) \leq c \left( \mathbb{E}\left( \sup_{[0,s]} |u^n|_{H^2}^2 \right) \right)^{1/2}.$$

On the other hand, it is easily seen from (6.3) that, for any $m \in \mathbb{N}$,

$$\mathbb{E}\left( \sup_{t \in [0,T]} |u^n(t)|_{H^{s_0}}^m \right) \leq C(m, T, R) \leq c$$

and we deduce

$$\mathbb{E}\left( \sup_{t \in [0,T]} A_1^n(t) \right) \leq c \left( \mathbb{E}\left( \sup_{t \in [0,T]} |u^n(t)|_{H^2}^2 \right) \right)^{1/2}.$$

It is now easy to gather all the preceding estimates and to conclude the proof of our claim.  □

STEP 2.   *For any $R \geq 1$, there exists a constant $C(R, T, \phi)$ such that, for any $n \in \mathbb{N}$,*

$$\mathbb{E}\left( \sup_{t \in [0,T]} |u_R^n(t)|_{\Sigma^2}^2 \right) \leq C(R, T, \phi).$$

PROOF. Thanks to the result of Step 1, it remains to estimate $V_2(u_R^n(t)) = \int_{\mathbb{R}^d} |x|^4 |u_R^n(s,x)|^2 \, dx$. We use the Itô formula, recalling as in Section 6.1 that the computations can be justified by a smoothing argument. As above, we omit to write the dependence on $R$. We obtain

$$dV_2(u^n) = 2\mathrm{Im} \int_{\mathbb{R}^d} |x|^4 \bar{u}^n \Delta u^n \, dx \, dt = -8\mathrm{Im} \int_{\mathbb{R}^d} |x|^2 \bar{u}^n x \cdot \nabla u^n \, dx \, dt.$$

We have, by Cauchy–Schwarz's inequality and integration by parts,

$$\left| \int_{\mathbb{R}^d} |x|^2 \bar{u}^n x \cdot \nabla u^n \, dx \right| \leq (V_2(u^n))^{1/2} |x \cdot \nabla u^n|$$

$$\leq c(V_2(u^n))^{1/2} ((V_2(u^n))^{1/4} |u^n|_{H^2}^{1/2} + |u^n|)$$

$$\leq \tfrac{1}{8} V_2(u^n) + c|u^n|_{H^2}^2.$$

We deduce

$$V_2(u^n(t)) \leq e^t \left[ V_2(u_0) + \int_0^t |u^n|_{H^2}^2 \, ds \right].$$

We conclude thanks to Step 1. $\quad\square$

STEP 3. *The laws of* $((u_R^n, W))_{n \in \mathbb{N}}$ *are tight in* $C([0,T]; \Sigma^\eta) \times C([0,T]; H_{\mathrm{loc}}^\eta(\mathbb{R}^d))$ *for any* $\eta < 2$.

PROOF. We use the integral form of the equation satisfied by $u_R^n$ to get an estimate on its modulus of continuity. Again, we omit to write the dependence on $R$ and all the constants appearing below do not depend on $n$. We have, for $t_1, t_2 \in [0,T]$ with $t_1 \leq t_2$,

$$S(-t_1)u^n(t_1) - S(-t_2)u^n(t_2)$$

$$= -i \int_{t_1}^{t_2} S(-s)[\theta(u^n)|u^n|^{2\sigma} u^n + f_n u^n + \Lambda(u^n)] \, ds$$

$$- i \int_{t_1}^{t_2} S(-s)\theta(u^n)u^n \, dW(s) + i \int_{t_1}^{t_2} S(-s)\theta(u^n)u^n \dot{W}_n \, ds.$$

Using the fact that $H^{s_0}(\mathbb{R}^d)$ is an algebra and that $S(-s)$ is an isometry, we deduce

$$\mathbb{E}(|S(-t_1)u^n(t_1) - S(-t_2)u^n(t_2)|_{H^{s_0}}^2)$$

$$\leq c|t_2 - t_1| + \mathbb{E}\left( \left| \int_{t_1}^{t_2} S(-s)\theta(u^n)u^n \dot{W}_n \, ds \right|_{H^{s_0}}^2 \right),$$

where again the constant $c$ is allowed to depend on $R$. The last term is decomposed as $A^n$ in the first step into the sum of a martingale term

$$\int_{t_1}^{t_2} S(-s_{k-1})\theta(u_S^n(s,s_{k-1}))u_S^n(s,s_{k-1})\dot{W}_n \, ds$$

and of a remaining term. It is not difficult to see that both have a second moment which can be majorized by a constant time $|t_2 - t_1|$. Therefore,

$$\mathbb{E}(|S(-t_1)u^n(t_1) - S(-t_2)u^n(t_2)|^2_{H^{s_0}}) \leq c|t_2 - t_1|.$$

Now, since $S(t)$ is an isometry and since, for $v$ in $H^2(\mathbb{R}^d)$, we have

$$(6.6) \qquad |(S(-t_1) - S(-t_2))v|_{H^{s_0}} \leq c|t_1 - t_2|^{2-s_0}|v|_{H^2},$$

as can be seen with the use of spatial Fourier transform and we easily obtain, thanks to Step 1,

$$\mathbb{E}(|u^n(t_1) - u^n(t_2)|^2_{H^{s_0}}) \leq c|t_2 - t_1|^{2(2-s_0)}.$$

Thus, the sequence $(u^n)_{n \in \mathbb{N}}$ is bounded in probability in $C([0,T]; \Sigma^2(\mathbb{R}^d)) \cap C^{2-s_0}([0,T]; H^{s_0}(\mathbb{R}^d))$. This space is compactly embedded into $C([0,T]; \Sigma^\eta(\mathbb{R}^d))$ for $\eta < 2$ and the tightness of $(u^n)_{n \in \mathbb{N}}$ follows.

We conclude thanks to the well-known identity

$$\mathbb{E}(|W(t_1) - W(t_2)|^2_{H^2}) = |t_1 - t_2| \, |\phi|^2_{\mathcal{L}^{0,2}_2}$$

and the compactness of the embedding of $H^2(\mathbb{R}^d)$ into $H^\eta_{\text{loc}}(\mathbb{R}^d)$. $\quad\square$

STEP 4. *Passage to the limit.*

Let us fix a pair of subsequences $(\varphi(n), \psi(n))_{n \in \mathbb{N}}$. We infer from Step 3 and the Prokhorov and Skorokhod theorems that there is a subsequence of $(u^{\varphi(n)}_R, u^{\psi(n)}_R, W)$ which we still denote by the same letters, a probability space $(\widetilde{\Omega}, \widetilde{\mathcal{F}}, \widetilde{\mathbb{P}})$ and random variables $(\tilde{u}^n_1, \tilde{u}^n_2, \widetilde{W}^n)$, $n \in \mathbb{N}$, $(\tilde{u}_1, \tilde{u}_2, \widetilde{W})$ with values in $C([0,T]; \Sigma^\eta) \times C([0,T]; \Sigma^\eta) \times C([0,T]; H^\eta_{\text{loc}}(\mathbb{R}^d))$ for any $\eta < 2$, such that for any $n \in \mathbb{N}$,

$$\mathcal{L}(\tilde{u}^n_1, \tilde{u}^n_2, \widetilde{W}^n) = \mathcal{L}(u^{\varphi(n)}_R, u^{\psi(n)}_R, W)$$

and such that

$$\tilde{u}^n_j \to \tilde{u}_j \qquad \text{as } n \to +\infty, \ \widetilde{\mathbb{P}} \text{ a.s. in } C([0,T]; \Sigma^\eta) \text{ for } \eta < 2, \ j = 1, 2,$$

$$\widetilde{W}^n \to \widetilde{W} \qquad \text{as } n \to +\infty, \ \widetilde{\mathbb{P}} \text{ a.s. in } C([0,T]; H^\eta_{\text{loc}}(\mathbb{R}^d)) \text{ for } \eta < 2.$$

Defining then

$$\widetilde{\mathcal{F}}_t = \sigma\{\tilde{u}_j(s), \widetilde{W}(s), \ 0 \leq s \leq t, \ j = 1, 2\}$$

and

$$\widetilde{\mathcal{F}}^n_t = \sigma\{\tilde{u}^n_j(s), \widetilde{W}^n(s), \ 0 \leq s \leq t, \ j = 1, 2\},$$

it is easily seen that $\widetilde{W}$ and $\widetilde{W}^n$ are Wiener processes associated, respectively, with $(\widetilde{\mathcal{F}}_t)_{t \geq 0}$ and $(\widetilde{\mathcal{F}}^n_t)_{t \geq 0}$, with covariance operator $\phi\phi^*$.

It can be checked that for $j = 1, 2$, we have

$$(6.7) \quad i d\tilde{u}_j^n - \Delta \tilde{u}_j^n \, dt - \theta_R(\tilde{u}_j^n)(|\tilde{u}_j^n|^{2\sigma} \tilde{u}_j^n + f_n \tilde{u}_j^n) \, dt$$

$$= \theta_R(\tilde{u}_j^n)\tilde{u}_j^n \, d\widetilde{W}^n + \Lambda_{R,\phi}(\tilde{u}_j^n) \, dt - \theta_R(\tilde{u}_j^n)\tilde{u}_j^n \dot{\widetilde{W}}_n \, dt$$

and

$$\tilde{u}_j^n(0) = u_0,$$

where $\dot{\widetilde{W}}_n$ is defined as $\dot{W}_n$ but replacing $W$ by $\widetilde{W}^n$. The difficulty now is to prove that the right-hand side of (6.7) goes to zero. In fact, we have the following result whose proof is very technical and is postponed to Step 6.

LEMMA 6.1.    *For $j = 1, 2$ and any $t \in [0, T]$,*

$$\int_0^t \theta_R(\tilde{u}_j^n)\tilde{u}_j^n \, d\widetilde{W}^n(s) + \int_0^t \Lambda_{R,\phi}(\tilde{u}_j^n) \, ds - \int_0^t \theta_R(\tilde{u}_j^n)\tilde{u}_j^n \dot{\widetilde{W}}_n \, ds \to 0,$$

*in $L^1(\tilde{\Omega}, H^{s_0}(\mathbb{R}^d))$.*

It is not difficult to let $n \to \infty$ in the left-hand side of (6.7) and, thanks to Lemma 6.1, we deduce that, for $j = 1, 2$,

$$\begin{cases} i d\tilde{u}_j - \Delta \tilde{u}_j \, dt - \theta_R(\tilde{u}_j)(|\tilde{u}_j|^{2\sigma} \tilde{u}_j + f\tilde{u}_j) \, dt = 0, \\ \tilde{u}_j(0) = u_0. \end{cases}$$

It can be seen that this equation has at most one solution so that the limits of $\tilde{u}_1^n$ and $\tilde{u}_2^n$ coincide:

$$\tilde{u}_1 = \tilde{u}_2.$$

We deduce that for each pair of subsequences $(\varphi(n), \psi(n))_{n \in \mathbb{N}}$, there is a subsequence of $(u_R^{\varphi(n)}, u_R^{\psi(n)})$ which converges in law to a random variable $(\tilde{u}_1, \tilde{u}_1)$. We now use the following elementary lemma, which was first used by Gyöngy and Krylov in [16]:

LEMMA 6.2.    *Let $Z_n$ be a sequence of random elements in a Polish space $E$ equipped with the Borel $\sigma$-algebra. Then $Z_n$ converges in probability to an $E$-valued random element if and only if for every pair of subsequences $(Z_{\varphi(n)}, Z_{\psi(n)})$, there is a subsequence of $(Z_{\varphi(n)}, Z_{\psi(n)})$ which converges in law to a random element supported on the diagonal $\{(x, y) \in E \times E, \ x = y\}$.*

It follows that the sequence $(u_R^n)_{n \in \mathbb{N}}$ converges in probability in $C([0, T], \Sigma^\eta)$ to a random variable $u_R$. Moreover, proceeding as above, we can let $n \to \infty$ in (6.3) and see that $u_R$ satisfies the equation

$$\begin{cases} i d u_R - \Delta u_R \, dt - \theta_R(u_R)(|u_R|^{2\sigma} u_R + f u_R) \, dt = 0, \\ \tilde{u}_R(0) = u_0. \end{cases}$$

Step 5.  *Conclusion.*

It is now important to emphasize the dependence on $R$. Let us define $\tau_R = \min\{t \in [0,T], \ |u_R(t)|_{H^{s_0}} \geq R\}$ and $\tau_R^n = \min\{t \in [0,T], \ |u_R^n(t)|_{H^{s_0}} \geq R\}$. On $[0,\tau_R]$, $u_R$ satisfies the deterministic nonlinear Schrödinger equation (5.2) whose unique solution is given by $U$. Thus, if $R_0 = |U|_{L^\infty(0,T_2;H^{s_0}(\mathbb{R}^d))}$ and $R > R_0$, then $\tau_R > T_2$ and $u_R = U$ on $[0,T_2]$. Similarly, on $[0,\tau_R^n]$, $u_R^n$ and $u^n$, the solution of (5.6) coincide. Moreover, $\tau_R^n$ and $\tau_R$ increase with $R$ and

$$\tau_R \to t^* \quad \text{and} \quad \tau_R^n \to \tau_n^*(u_0) \qquad \text{a.s.}$$

when $R \to \infty$, where $t^*$ is the existence time of $U$. It follows, since $t^* > T_2$, that

$$\mathbb{P}(\tau_n^*(u_0) \leq T_2) \to 0$$

as $n \to \infty$. Hence, (5.8) holds. Moreover, for any $\epsilon \leq 1$,

$$\mathbb{P}\bigg(\tau_n^*(u_0) > T_2 \text{ and } \sup_{t \in [0,T_2]} |u^n - U|_\Sigma > \epsilon\bigg)$$

$$\leq \mathbb{P}\bigg(\tau_n^*(u_0) > T_2 \text{ and } \sup_{t \in [0,T_2 \wedge \tau_{R_0+1}^n]} |u^n - U|_\Sigma > \epsilon\bigg)$$

$$\leq \mathbb{P}\bigg(\sup_{t \in [0,T]} |u_{R_0+1}^n(t) - u_{R_0+1}(t)|_\Sigma > \epsilon\bigg),$$

where we have used $\tau_{R_0+1} > T_2$. Now, this last term goes to 0 as $n$ goes to infinity by Step 4. We deduce that (5.9) holds.

Step 6.  *Proof of Lemma 6.1.*

In order to lighten the notation, we omit writing the dependence on $R$. Also we write $\Lambda$ instead of $\Lambda_{R,\phi}$, $u^n$ instead of $\tilde{u}_j^n$ and $W^n$ instead of $\widetilde{W}^n$. Thus, we have to prove

$$I_n(t) = \int_0^t \theta(u^n)u^n \, dW^n(s) + \int_0^t \Lambda(u^n) \, ds - \int_0^t \theta(u^n)u^n \dot{W}_n \, ds \to 0,$$

when $n \to \infty$ for any $t \in [0,T]$ in $L^1(\widetilde{\Omega}, H^{s_0}(\mathbb{R}^d))$. We split $I_n$ as the sum of $I_n^1$ and $I_n^2$, where

$$I_n^1(t) = \sum_{k=0}^{N-1} \int_{t_k}^{t_{k+1}} \theta(u^n(s))u^n(s) \, dW^n(s) - \int_{t_k}^{t_{k+1}} \theta(u_S^n(s,t_{k-1}))u_S^n(s,t_{k-1})\dot{W}_n \, ds$$

and

$$I_n^2(t) = \sum_{k=0}^{N-1} \int_{t_k}^{t_{k+1}} (\theta(u_S^n(s, t_{k-1})) u_S^n(s, t_{k-1}) - \theta(u^n(s)) u^n(s)) \dot{W}_n \, ds$$
$$+ \int_{t_k}^{t_{k+1}} \Lambda(u^n(s)) \, ds.$$

We have set $t_k = k\Delta t$, $k = 0, \ldots, N-1$, $t_N = T$, where $N$ is the integer part of $t/\Delta t$, and $t_{-1} = 0$. We have again used the notation $u_S^n(s, t_{k-1}) = S(s - t_{k-1}) u^n(t_{k-1})$. Note that, for $k = 0, \ldots, N-1$,

$$\int_{t_k}^{t_{k+1}} \theta(u_S^n(s, t_{k-1})) u_S^n(s, t_{k-1}) \dot{W}_n \, ds$$
$$= \int_{t_{k-1}}^{t_k} \frac{1}{\Delta t} \left( \int_{t_k}^{t_{k+1}} \theta(u_S^n(\tau, t_{k-1})) u_S^n(\tau, t_{k-1}) \, d\tau \right) d(\phi P_n W_c^n)(s),$$

where $W_c^n = \phi^{-1} W^n$ is a cylindrical Wiener process. It follows

$$I_n^1(t) = \sum_{k=0}^{N-1} \int_{t_{k-1}}^{t_k} \frac{1}{\Delta t} \left( \int_{t_k}^{t_{k+1}} (\theta(u^n(s)) u^n(s) \right.$$
$$\left. - \theta(u_S(\tau, t_{k-1})) u_S(\tau, t_{k-1})) \, d\tau \right) dW^n(s)$$
$$+ \int_{t_{N-1}}^{t_N} \theta(u^n(s)) u^n(s) \, dW^n(s)$$
$$+ \sum_{k=0}^{N-1} \int_{t_{k-1}}^{t_k} \theta(u^n(s)) u^n(s) \, d(\phi(I - P_n) W_c^n).$$

Then, using similar arguments as in Step 1,

$$\mathbb{E}(|I_n^1(t)|_{H^{s_0}}^2)$$
$$\leq \sum_{k=0}^{N-1} \mathbb{E}\left( \int_{t_{k-1}}^{t_k} \left| \frac{1}{\Delta t} \left( \int_{t_k}^{t_{k+1}} (\theta(u^n(s)) u^n(s) \right. \right. \right.$$
$$\left. \left. \left. - \theta(u_S(\tau, t_{k-1})) u_S(\tau, t_{k-1})) \, d\tau \right) \phi \right|_{\mathcal{L}_2^{0,s_0}}^2 \, ds \right)$$

(6.8)
$$+ \int_0^t \mathbb{E}(|\theta(u^n(s)) u^n(s) \phi(I - P_n)|_{\mathcal{L}_2^{0,s_0}}^2) \, ds$$
$$+ \int_{t_{N-1}}^{t_N} \mathbb{E}(|\theta(u^n(s)) u^n(s) \phi|_{\mathcal{L}_2^{0,s_0}}^2) \, ds$$
$$\leq c(\Delta t) + c|\phi(I - P_n)|_{\mathcal{L}_2^{0,s_0}}$$

and this latter term goes to 0 as $n \to \infty$, since $\phi \in \mathcal{L}_2^{0,s_0}$.

Now, estimate (6.5) is not sufficient to handle $I_n^2(t)$. Let us define

$$R(s, t_{k-1}) = -i \int_{t_{k-1}}^s S(s-\tau)[\theta(u^n)(|u^n|^{2\sigma}u^n + f_n u^n) - \Lambda(u^n)] \, d\tau$$

$$- i \int_{t_{k-1}}^s S(s-\tau)[\theta(u^n)u^n \, dW^n(\tau)]$$

$$+ i \int_{t_{k-1}}^s \theta(u_S^n(s, t_{k-1}))u_S^n(s, t_{k-1}) \, dW^n(\tau)$$

$$+ i \int_{t_{k-1}}^s S(s-\tau)[\theta(u^n)u^n \dot{W}_n \, d\tau]$$

$$- i \int_{t_{k-1}}^s \theta(u_S^n(s, t_{k-1}))u_S^n(s, t_{k-1}) \dot{W}_n \, d\tau$$

so that

$$u^n(s) - u_S^n(s, t_{k-1})$$

$$= R(s, t_{k-1}) - i\theta(u_S^n(s, t_{k-1}))u_S^n(s, t_{k-1}) \int_{t_{k-1}}^s dW^n(\tau)$$

$$+ i\theta(u_S^n(s, t_{k-1}))u_S^n(s, t_{k-1}) \int_{t_{k-1}}^s \dot{W}_n \, d\tau.$$

We can prove, using Cauchy–Schwarz's inequality, Step 1, (6.5) and (6.6), that for any $\delta < 2$,

$$\mathbb{E}(|R(s, t_{k-1})|_{H^{s_0}}^\delta) \leq c\Delta t^{\delta(2-s_0+1/2)}.$$

It follows

$$(6.9) \qquad \mathbb{E}(|R(s, t_{k-1})|_{H^{s_0}}|\dot{W}_n|_{H^{s_0}}) \leq c\Delta t^{2-s_0}.$$

Next, we rewrite $I_n^2$ as

$$I_n^2(t) = \sum_{k=0}^{N-1} \int_{t_k}^{t_{k+1}} (\theta(u_S^n(s, t_{k-1}))u_S^n(t_{k-1}) - \theta(u^n(s))u^n(s)$$

$$- \sigma_\theta'(u_S^n(s, t_{k-1})) \cdot (u_S^n(s, t_{k-1}) - u^n(s)))\dot{W}_n \, ds$$

$$- \sum_{k=0}^{N-1} \int_{t_k}^{t_{k+1}} \sigma_\theta'(u_S^n(s, t_{k-1})) \cdot R(s, t_{k-1})\dot{W}_n \, ds$$

$$+ \sum_{k=0}^{N-1} \int_{t_k}^{t_{k+1}} (\Lambda(u^n(s)) - \Lambda(u_S^n(s, t_{k-1}))) \, ds$$

$$+ \sum_{k=0}^{N-1} \int_{t_k}^{t_{k+1}} \theta(u_S^n(s, t_{k-1})) \sigma_\theta'(u_S^n(s, t_{k-1}))$$

$$\cdot \left( i u_S^n(s, t_{k-1}) \int_{t_{k-1}}^{s} dW^n(\tau) \right) \dot{W}_n \, ds$$

$$- \sum_{k=0}^{N-1} \int_{t_k}^{t_{k+1}} \theta(u_S^n(s, t_{k-1})) \sigma_\theta'(u_S^n(s, t_{k-1}))$$

$$\cdot \left( i u_S^n(s, t_{k-1}) \int_{t_{k-1}}^{s} \dot{W}_n \, d\tau \right) \dot{W}_n \, ds$$

$$+ \sum_{k=0}^{N-1} \int_{t_k}^{t_{k+1}} \Lambda(u_S^n(s, t_{k-1})) \, ds,$$

where $\sigma_\theta'(v) \cdot w$ is the differential in the direction $w$ of the mapping $v \mapsto \theta(v)v$ which is $C^\infty$ from $H^{s_0}$ into itself with all derivatives bounded. By the Taylor formula and (6.5), the first term clearly goes to zero in $L^1(\Omega; H^{s_0})$ as $(\Delta t)^{1/2}$. This is also the case for the third term. For the second term, we use (6.9) and see that it goes to zero in $L^1(\Omega; H^{s_0})$ as $(\Delta t)^{2-s_0}$.

Let us denote by $J(t)$ the sum of the remaining terms. Using that

$$\Lambda(u) = \tfrac{1}{2} \theta(u) \sum_{\ell \in \mathbb{N}} \sigma_\theta'(u) \cdot (-i u \phi e_\ell) \phi e_\ell,$$

we get

$$J(t) = \sum_{k=0}^{N-1} \int_{t_k}^{t_{k+1}} \theta(u_S^n(s, t_{k-1}))$$

$$\times \left[ \sigma_\theta'(u_S^n(s, t_{k-1})) \right.$$

$$\cdot \left( i u_S^n(s, t_{k-1}) \left( \int_{t_{k-1}}^{s} dW^n(\tau) - \int_{t_{k-1}}^{s} \dot{W}_n \, d\tau \right) \right) \dot{W}_n$$

$$\left. + \tfrac{1}{2} \sum_{\ell \in \mathbb{N}} \sigma_\theta'(u_S^n(s, t_{k-1})) \cdot (-i u_S^n(s, t_{k-1}) \phi e_\ell) \phi e_\ell \right] ds.$$

Now, recalling the definition of $\theta$, we have

$$\sigma_\theta'(v) \cdot w = \frac{2}{R^2} \vartheta' \left( \frac{|v|_{H^{s_0}}^2}{R^2} \right) (v, w)_{H^{s_0}} v + \vartheta \left( \frac{|v|_{H^{s_0}}^2}{R^2} \right) w,$$

and we can write

$$J(t) = J_1(t) + J_2(t),$$

with

$$J_1(t) = \sum_{k=0}^{N-1} \int_{t_k}^{t_{k+1}} i\theta^2(u_S^n(s, t_{k-1})) u_S^n(s, t_{k-1})$$

$$\times \left[ \left( \int_{t_{k-1}}^s dW^n(\tau) - \int_{t_{k-1}}^s \dot{W}_n \, d\tau \right) \dot{W}_n - \tfrac{1}{2} F_\phi \right] ds$$

and

$$J_2(t) = \frac{2}{R^2} \sum_{k=0}^{N-1} \int_{t_k}^{t_{k+1}} \theta(u_S^n(s, t_{k-1})) \vartheta' \left( \frac{|u_S^n(s, t_{k-1})|_{H^{s_0}}^2}{R^2} \right)$$

$$\times \left[ \left( u_S^n(s, t_{k-1}), i u_S^n(s, t_{k-1}) \right. \right.$$

$$\times \left. \left( \int_{t_{k-1}}^s dW^n(\tau) - \int_{t_{k-1}}^s \dot{W}_n \, d\tau \right) \right)_{H^{s_0}} \dot{W}_n$$

$$- \frac{1}{2} \sum_{\ell \in \mathbb{N}} (u_S^n(s, t_{k-1}), i u_S^n(s, t_{k-1}) \phi e_\ell)_{H^{s_0}} \phi e_\ell \right]$$

$$\times u_S^n(s, t_{k-1}) \, ds.$$

Writing

$$\left( \int_{t_{k-1}}^s dW^n(\tau) - \int_{t_{k-1}}^s \dot{W}_n \, d\tau \right) \dot{W}_n(s) - \tfrac{1}{2} F_\phi$$

$$= \left[ \left( \int_{t_{k-1}}^s dW^n(\tau) \right) \dot{W}_n(s) - F_\phi \right] - \left[ \left( \int_{t_{k-1}}^s \dot{W}_n(\tau) \, d\tau \right) \dot{W}_n(s) - \tfrac{1}{2} F_\phi \right],$$

it is then the result of tedious—but standard—computations, based on the independence of the increments of the Wiener process, that each bracket in the right-hand side above gives, when inserted into $J_1(t)$, a term which goes to 0 in $L^2(\Omega; L^2(\mathbb{R}^d))$, hence, also in $L^2(\Omega; H^{s_0}(\mathbb{R}^d))$—since it is clearly bounded in $L^2(\Omega; H^2(\mathbb{R}^d))$.

It can be proved in the same way that $J_2(t)$ goes to 0 in $L^2(\Omega; H^{s_0}(\mathbb{R}^d))$ as $\Delta t$ goes to 0 (or $n$ goes to infinity).

Therefore, $I_n^2(t)$ goes to zero in $L^1(\Omega; H^{s_0})$ and, recalling (6.8), Lemma 6.1 is proved.

CNRS et Université Paris-Sud
UMR 8628
Bâtiment 425
Université Paris-Sud
91405 Orsay Cedex
France
e-mail: anne.debouard@math.u-psud.fr

ENS de Cachan
Antenne de Bretagne
Campus de Ker Lann
Avenue R. Schuman
35170 Bruz Cedex
France
e-mail: arnaud.debussche@bretagne.ens-cachan.fr